\documentclass[a4paper,review]{cas-sc}

\usepackage[section]{placeins} 
\usepackage{amssymb}
\usepackage{amsthm}
\usepackage{amsmath}
\usepackage{upgreek}
\usepackage{pdflscape}
\usepackage{listings}
\usepackage{multirow}
\usepackage[percent]{overpic}
\usepackage{multicol}
\usepackage{color}
\usepackage{stmaryrd}
\usepackage{xfrac}
\usepackage[capitalise]{cleveref}
\usepackage{booktabs}
\usepackage[section]{placeins} 
\usepackage[section]{algorithm}
\usepackage{calc}
\usepackage[titletoc]{appendix}
\usepackage{siunitx}
\usepackage[sort&compress,square,numbers]{natbib}

\usepackage{graphicx}
\usepackage{graphics}
\usepackage{wrapfig}
\usepackage{float}
\usepackage{subfig}
\usepackage[percent]{overpic}
\usepackage{varwidth}
\usepackage{tikz}
\usepackage{pgfplots}
\usepackage{adjustbox}
\usetikzlibrary{arrows,matrix,positioning,fit}
\usetikzlibrary{shapes,positioning}
\usetikzlibrary{backgrounds}
\usepackage{tikz-layers}
\usepgfplotslibrary{fillbetween}
\usetikzlibrary{intersections}
\pgfplotsset{compat=1.14}
\usepgfplotslibrary{colorbrewer}
\usepgfplotslibrary{patchplots}
\usepgfplotslibrary[colorbrewer]
\usetikzlibrary{pgfplots.colorbrewer}
\usetikzlibrary[pgfplots.colorbrewer]
\usepgfplotslibrary{units}
\usetikzlibrary{spy}
\usepackage{pgfplotstable}
\usepackage{arrayjobx}
\graphicspath{ {./figs/} }
\usetikzlibrary{external}

\newcommand{\multiclip}[3]{
\newarray\temp%
\readarray{temp}{#2}%
\foreach \x in {1,...,#1}%
{ \begin{scope}%
        \pgfmathtruncatemacro{\xt}{\x}%
        \temp(\xt)%
        #3%
    \end{scope}%
}%
\delarray\temp%
}

\newlength\myheight
\newlength\mydepth
\settototalheight\myheight{Xygp}
\newcommand*\inlinegraphics[1]{%
  \settototalheight\myheight{Xygp}%
  \settodepth\mydepth{Xygp}%
  \raisebox{-\mydepth}{\includegraphics[height=\myheight]{#1}}%
}
\newcommand\orcid[1]{\href{https://orcid.org/#1}{\inlinegraphics{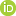}}}

\makeatletter
\def\BState{\State\hskip-\ALG@thistlm}
\makeatother

\newtheorem{theorem}{Theorem}[section]
\newtheorem{corollary}{Corollary}[theorem]
\newtheorem{lemma}[theorem]{Lemma}
\newdefinition{definition}{Definition}[section]

\newtheorem{example}{Example}[section]

\newcommand{\etal}{et~al.}

\newcommand\ob[1]{\overline{\mathbf{#1}}}

\newcommand\px[2]{\frac{\partial #1}{\partial {#2}}}

\newcommand\dx[2]{\frac{\mathrm{d} #1}{\mathrm{d} #2}}

\newcommand\pxvar[2]{\partial_{#2} #1}

\newcommand{\half}{{\frac{1}{2}}}
\newcommand{\shalf}{{\sfrac{1}{2}}}

\begin{document}

\title[mode=title]{A Riemann Difference Scheme for Shock Capturing in Discontinuous Finite Element Methods}
\shorttitle{A Riemann Difference Scheme for Shock Capturing in Discontinuous Finite Element Methods}
\shortauthors{T. Dzanic~\etal}

\author[1]{T. Dzanic}[orcid=0000-0003-3791-1134]
\cormark[1]
\cortext[cor1]{Corresponding author}
\ead{tdzanic@tamu.edu}
\author[1]{W. Trojak}[orcid=0000-0002-4407-8956]
\author[1]{F. D. Witherden}[orcid=0000-0003-2343-412X]

\address[1]{Department of Ocean Engineering, Texas A\&M University, College Station, TX 77843}

\begin{abstract}
	We present a novel structure-preserving numerical scheme for discontinuous finite element approximations of nonlinear hyperbolic systems. The method can be understood as a generalization of the Lax--Friedrichs flux to a high-order staggered grid and does not depend on any tunable parameters. Under a presented set of conditions, we show that the method is conservative and invariant domain preserving. Numerical experiments on the Euler equations show the ability of the scheme to resolve discontinuities without introducing excessive spurious oscillations or dissipation. 
\end{abstract}

\begin{highlights}
\item A novel staggered grid spectral element method is presented for conservation laws 
\item The method is proved to be invariant domain preserving under set conditions
\item Numerical test cases with discontinuities show the high resolution and low dissipation of the method
\end{highlights}




\maketitle

\section{Introduction}
\label{sec:intro}

Hyperbolic conservation laws in continuum dynamics govern the behavior of many systems of interest to scientists and engineers. It has been known for considerable time that, when nonlinear, these systems can produce discontinuities in finite time even with smooth initial conditions \citep{Hopf1950}. For the development of numerical schemes under these circumstances, Godunov's theorem, \citep{Godunov1959}, presents the crux of the matter -- linear schemes for the monotonic solution of hyperbolic equations can be at most first-order accurate. Consequently, the formation of high-resolution approximations for solutions of nonlinear hyperbolic systems of equations poses a significant challenge.

Discontinuous finite element methods (FEM) have grown in prevalence with the increased adoption of highly parallel computation hardware; however, the issues presented by discontinuous solutions is one of the limitations preventing the widespread adoption of these methods for the industry. Although these methods provide a mathematically robust manner in which to achieve an arbitrarily high order of accuracy, when confronted by discontinuous solutions or large gradients, the presence of Gibbs phenomena can cause the solution to become aphysical or diverge \citep{Lax2006}. Several strategies have been proposed to permit the use of FEM on solutions that exhibit discontinuities. The classical and most ubiquitous methods are variations on the artificial viscosity (AV) approach, first proposed by \citet{VonNeumann1950} and subsequently modified to the spectral vanishing viscosity method by \citet{Tadmor1990}. Originally intended for finite difference and pseudo-spectral methods, they have subsequently been applied to spectral difference and discontinuous Galerkin methods \citep{Glaubitz2017,Persson2006}. The issue presented to practitioners by these AV methods is the use of tunable parameters controlling the amount of dissipation added.

A comparatively recent development has been the invariant domain preserving graph viscosity method of \citet{Guermond2016} which was later generalized to an abstract numerical setting by \citet{Guermond2019}. This technique defines a low-order approximation within a high-order framework that has the unique property of preserving all convex invariants of the hyperbolic system in question \citep{Glimm1965,Chueh1977,Hoff1985}. This may be thought of as a more physical interpretation of the total variation diminishing (TVD) property and is fundamentally linked to the entropy condition for solutions of hyperbolic systems. Invariant domain preserving graph viscosity has provable properties which are of interest to those seeking to highly resolve discontinuous solutions and does so without the need for tunable parameters, but it comes at the expense of being highly dissipative. The results may be considerably improved when combined with a convex limiting procedure---similar to the flux-corrected transport~(FCT) method of \citet{Boris1997}---and also by applying the entropy viscosity methods of \citet{Guermond2011}. However, this introduces a similar issue to that confronted by many schemes for resolving discontinuities: algorithmic complexity. The aim of this work then is to define a scheme that attempts to be algorithmically simple, devoid of parameterization, and posses lower dissipation properties than other approaches. 

Using the idea of the invariant domain preserving methods, we propose a novel scheme with the aim of significantly reducing the overall dissipation such that the resulting scheme is suitable for scale resolving simulations. This method builds off of the staggered grid approach of \citet{Kopriva1996} and utilizes approximate solutions of the Riemann problem to provide a physically proportionate amount of diffusion. In \cref{sec:prelim}, we outline the general hyperbolic system and some pertinent mathematical findings of previous works as a preliminary to our main results. The proposed scheme and the main theorems of this paper are presented in \cref{sec:scheme} with implementation details given in \cref{sec:numerical}. The scheme is investigated through numerical experiments in \cref{sec:results}, and conclusions are drawn in \cref{sec:conclusions}.
\section{Preliminaries}\label{sec:prelim}
 This work pertains to the solution of hyperbolic conservation laws of the form
\begin{equation}\label{eq:gen_hype}
     \begin{cases}
        \pxvar{\mathbf{u}}{t} + \mathbf{\nabla}\cdot\mathbf{F}(\mathbf{u}) = 0, \quad \mathrm{for}\ (\mathbf{x}, t)\in\mathcal{D}\times\mathbb{R}_+, \\
        \mathbf{u}(\mathbf{x}, 0)=\mathbf{u}_0, \quad \mathrm{for}\ \mathbf{x}\in\mathcal{D},
    \end{cases}
\end{equation}
where the solution $\mathbf{u}\in\mathbb{R}^m$, the flux $\mathbf{F}(\mathbf{u})\in(\mathbb{R}^m)^d$, and $\mathcal{D}\subseteq\mathbb{R}^d$ for some arbitrary space dimension $d$. To reduce the complexity of the analysis, we impose periodic boundary conditions on $\mathcal{D}$. To aid in proving the numerical properties of the scheme, we will present some properties of convex sets, invariant sets and domains, and summation-by-parts.

\subsection{Convex Sets}
    \begin{definition}[Minimum Distance]\label{def:b_dist}
        For a closed set $X$, the shortest distance from a state $\mathbf{u}\in X$ to the boundary of the set, $\partial X$, is defined as 
        \begin{equation}
            G_X(\mathbf{u}) = \inf_{\mathbf{x}\in\partial X}\|\mathbf{u} - \mathbf{x}\|_X,
        \end{equation}
        for some norm $\|\cdot\|_X$ on $X$.
    \end{definition}

    \begin{lemma}[Combined Distance]\label{lem:convex_dist}
        For a closed convex subset $X$ of a finite vector space, $A$, and two states $\mathbf{u}_1, \mathbf{u}_2 \in X$, 
        \begin{equation}\label{eq:com_dist}
            G_X(a_1 \mathbf{u}_1 + a_2 \mathbf{u}_2) \geqslant \mathrm{min} \left [ G_X(\mathbf{u}_1), G_X(\mathbf{u}_2) \right]
        \end{equation}
        for any convex combination such that $a_1,a_2\geqslant0$ and $a_1+a_2=1$.
    \end{lemma}

    \begin{theorem}[Prolongation Factor]\label{thm:con_factor}
        For a closed convex subset $X$ of a finite vector space, $A$, and a vector $\mathbf{u}\in X$, there exists a prolongation factor $\beta$ defined as 
        \begin{equation}
            \beta = 1 \pm \frac{G_X(\mathbf{u})}{\|\mathbf{u}\|_X},
        \end{equation}
        such that $\beta\mathbf{u}\in X$.
    \end{theorem}
    \newproof{pot_cf}{Proof of \cref{thm:con_factor}}
    \begin{pot_cf}
        Consider the closed ball of radius $r$ centered at $\mathbf{u}$, $\mathbb{B}(\mathbf{u},r)=\{\mathbf{x}\in A: \|\mathbf{u}-\mathbf{x}\|_X\leqslant r\}$. Given that $\mathbf{u}\in X$, then, by construction, $\mathbb{B}$ is convex subset of $X$ iff $r\leqslant G_X(\mathbf{u})$. It follows that $\mathbf{u} + \mathbf{e}G_X(\mathbf{u})\in X$ for any $\|\mathbf{e}\|_X=1$, and therefore
        \begin{equation*}
            \mathbf{u} \pm \frac{G_X(\mathbf{u})\mathbf{u}}{\|\mathbf{u}\|_X} = \left(1 \pm \frac{G_X(\mathbf{u})}{\|\mathbf{u}\|_X} \right)\mathbf{u} = \beta\mathbf{u} \in X.
        \end{equation*}
        \qed
    \end{pot_cf}
    This fundamental theorem leads us to the following corollary for non-convex combinations of members within and outside of convex sets.
    \begin{corollary}[Convex and Non-convex Combinations]\label{cor:noninvariant_sum}
        For a closed convex subset $X$ of a finite vector space, $A$, a state $\mathbf{u}_1 \in X$, and any state $\mathbf{u}_2$, by \cref{thm:con_factor}, 
        \begin{equation}\label{eq:nonconvex_c_sum}
            \mathbf{u}_3 = \mathbf{u}_1 + \delta s\ \mathbf{u}_2 \in X \quad \mathrm{for} \quad |\delta s| \leqslant \frac{G_X(\mathbf{u}_1)}{\|\mathbf{u}_2\|_X}.
        \end{equation}
    \end{corollary}
    
\subsection{Invariant Sets}\label{ssec:invariant_sets}
    In this subsection, we will restate the Riemann problem and discuss its connections to invariant sets in the context of \cref{eq:gen_hype}. These preliminaries will take a similar format to those of \citet{Guermond2019} and draws on several other works \citep{Glimm1965,Chueh1977,Hoff1985}.
    
    The Riemann problem, expressed as  
    \begin{equation}\label{eq:riemann}
        \pxvar{\mathbf{u}}{t} + \pxvar{(\mathbf{F}(\mathbf{u})\cdot\mathbf{n})}{x} = 0, \quad  (x,t) \in \mathbb{R}\times\mathbb{R}_+, \quad \mathbf{u}(x,0) = \begin{cases}
        \mathbf{u}_l, &\mbox{if } x<0\\
        \mathbf{u}_r, &\mbox{if } x>0
        \end{cases}, 
    \end{equation}
    for some normal vector $\mathbf{n}$, can be seen as a restriction to the broad class of equations defined by \cref{eq:gen_hype}. 
    From \citet{Dafermos2010_9} and \citet{Lax1957}, there is a unique self-similar solution $\mathbf{v}(\mathbf{u}_l,\mathbf{u}_r, \mathbf{n}, x/t)$ to this problem for a genuinely hyperbolic system with sufficiently small $|\mathbf{u}_l - \mathbf{u}_r|$. We state that there is an \emph{admissible} set, $\mathcal{A} \in \mathbb{R}^d$, such that the solution $\mathbf{v}(\mathbf{u}_l,\mathbf{u}_r, \mathbf{n}, x/t)\in\mathcal{A}$ for any $(\mathbf{u}_l,\mathbf{u}_r)\in\mathcal{A}$. If this system has a maximum absolute wavespeed, denoted by $\lambda_\mathrm{max}$, then the solution is $\mathbf{u}=\mathbf{u}_l$ for $x/t<-\lambda_\mathrm{max}$ and $\mathbf{u}=\mathbf{u}_r$ for $x/t>\lambda_\mathrm{max}$.
    
    \begin{lemma}[Average Riemann Solution] \label{rieman_lemma}
        Let the average Riemann solution over the Riemann fan, $\ob{v}$, for the Riemann problem in \cref{eq:riemann} be 
        \begin{equation}
            \ob{v}(\mathbf{u}_l,\mathbf{u}_r,\mathbf{n}, t) = \int^\shalf_{-\shalf}\mathbf{v}(\mathbf{u}_l,\mathbf{u}_r,\mathbf{n}, x/t)\ \mathrm{d}x.
        \end{equation}
        For a sufficiently small $t$ such that $\lambda_\mathrm{max}t \leqslant \frac{1}{2}$,
        \begin{subequations}
            \begin{align}
                \ob{v} &= \half(\mathbf{u}_r + \mathbf{u}_l) - t\left(\mathbf{F}(\mathbf{u}_r) - \mathbf{F}(\mathbf{u}_l)\right)\cdot\mathbf{n},\\
                \sigma(\ob{v}) &\leqslant \half(\sigma(\mathbf{u}_r) + \sigma(\mathbf{u}_l)) - t\left(\mathbf{\Sigma}(\mathbf{u}_r) - \mathbf{\Sigma}(\mathbf{u}_l)\right)\cdot\mathbf{n},
            \end{align}
        \end{subequations}
        for an entropy-flux pair $(\sigma, \mathbf{\Sigma})$. See \citet{Dafermos2008}.
    \end{lemma}
    
    \begin{definition}[Invariant Sets]\label{def:invariant_set}
        A set $\mathcal{B}\subset\mathcal{A}$ is said to be invariant w.r.t. \cref{eq:gen_hype} if for any pair of initial conditions $(\mathbf{u}_l,\mathbf{u}_r) \in \mathcal{B}\times\mathcal{B}$ and any unit vector $\mathbf{n}\in\mathbb{B}^{d-1}(\mathbf{0},1)$,  $\ob{v}$ remains in $\mathcal{B}$ for any $t > 0$ such that $\lambda_\mathrm{max}t \leqslant \frac{1}{2}$. 
    \end{definition}
    
    From the work of \citet{Hoff1985}, we see that for genuinely nonlinear equations, the invariant set is a convex set. We invoke this assumption as it allows for the use of the previously presented properties of convex sets. A further important notion to define is that of an invariant domain for which we take the same definition as \citet{Guermond2019}. 
    
    \begin{definition}[Invariant Domain]
        Let $\mathbf{U}=(\mathbf{u}_1,\dots,\mathbf{u}_N)$ for some positive integer $N$. A convex invariant set $\mathcal{B}$ is said to be an invariant domain for a mapping $R:(\mathbb{R}^d)^N\rightarrow(\mathbb{R}^d)^N$ iff $R(\mathbf{U})\in\mathcal{B}^N$ for any $\mathbf{U}\in\mathcal{B}^N$.
    \end{definition}
    
    \begin{example}[Compressible Euler Equations]\label{ex:euler}
        The compressible Euler equations may be written in the form of \cref{eq:gen_hype} as
        \begin{equation}
            \mathbf{u} = \begin{bmatrix}
                    \rho \\ \boldsymbol{\rho v} \\ E
                \end{bmatrix}, \quad  \mathbf{F} = \begin{bmatrix}
                    \boldsymbol{\rho v}\\
                    \boldsymbol{\rho v}\otimes\mathbf{v} + p\mathbf{I}\\
                (E+p)\mathbf{v}
            \end{bmatrix},
        \end{equation}
        where $\rho$ is the density, $\boldsymbol{\rho v}$ is the momentum, $E$ is the total energy, $p = (\gamma-1)\left(E - \half\rho\|\mathbf{v}\|_2^2\right)$ is the pressure, and $\gamma$ is the ratio of specific heat capacities. The symbol $\mathbf{I}$ denotes the identity matrix in $\mathbb{R}^{d\times d}$ and $\mathbf{v} = \boldsymbol{\rho v}/\rho$ denotes the velocity. From \citet{Guermond2016}, for a specific internal energy $e(\mathbf{u})  := E/\rho - \frac{1}{2}\|\mathbf{v}\|_2^2$ and specific physical entropy $s(\mathbf{u})$ such that $-s(e, \rho^{-1})$ is a strictly convex function, the set 
        \begin{equation}
            \mathcal{A} := \{(\rho, \boldsymbol{\rho v}, E)\ |\ \rho \geqslant 0, e(\mathbf{u})  \geqslant 0, s(\mathbf{u}) \geqslant s_0\}
        \end{equation}
        in an invariant set for the Euler system for any $s_0 \in \mathbb{R}$. An example of this invariant set is shown in \cref{fig:convex_set}.
        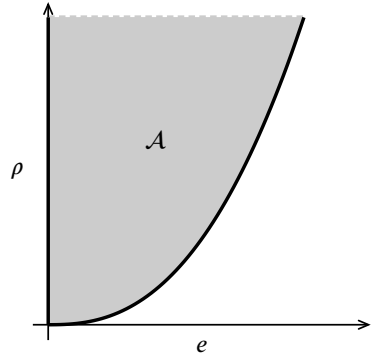
\begin{figure}[tbhp]
            \centering
            \adjustbox{width=0.3\linewidth,valign=b}{    \begin{tikzpicture}[scale = 0.5]

		\draw[name path=B,black!00] plot[smooth,samples=32,domain=0:8.325532074] (\x,{0.05*(\x)^2.5)});
        \draw[name path=A,black!00] (0,0) -- (0,10) -- (8.325523074,10);
        \draw[name path=A,ultra thick,densely dashed,gray!40] (0,10.01) -- (8.325523074,10.01);
        \draw[black,thick,->,>=angle 45] (0,-0.5) -- (0,10.5) node[midway,left,xshift=-2ex]{\large $\rho$};
        \draw[black,thick,->,>=angle 45] (-0.5,0) -- (10.5,0) node[midway,below,yshift=-1ex]{\large $e$};
        \tikzfillbetween[of=A and B]{gray!40};
        
        \node[black] at (3.5,6) {\large $\mathcal{A}$};
    
		\draw[ultra thick,black] plot[smooth,samples=32,domain=0:8.325532074] (\x,{0.05*(\x)^2.5)});
        \draw[black,ultra thick] (0,0) -- (0,10);
	\end{tikzpicture}}
            \caption{\label{fig:convex_set} Example of the invariant set $\mathcal{A}$ for the Euler system for $s(\mathbf{u}) = p \rho^{-\gamma}$, $s_0 > 0$, and $\gamma = 1.4$.}
        \end{figure}
    \end{example}
    

\subsection{Summation-By-Parts}
    For an $n$ dimensional discretization $\mathbf{x}$ of the reference space $\bar{\Omega}=[-1,1]$, let $\mathbf{v},\mathbf{w}\in\mathbb{R}^n$ be such that $\mathbf{v}=v(\mathbf{x})$, $\mathbf{w}=w(\mathbf{x})$ for some functions $v,w\in H^1$. The mass matrix, $\mathbf{M}\in\mathbb{R}^{n\times n}$, is defined such that
    \begin{equation}
        \mathbf{v}^T\mathbf{M}\mathbf{w} \approx  \langle v, w\rangle,
    \end{equation}
    where $\langle \cdot ,\cdot\rangle$ denotes the inner product. A discrete differentiation operator, $\mathbf{D}\in\mathbb{R}^{n\times n}$, may also be defined as
    \begin{equation}
        \mathbf{D}\mathbf{v} \approx \pxvar{v}{x}|_\mathbf{x}.
    \end{equation}
    Furthermore, we require a boundary projection operator $\mathbf{P}$ such that $\mathbf{P}\mathbf{v}=[v(-1),v(1)]^T$ and the boundary operator $\mathbf{B}=\mathrm{diag}{(-1,1)}$. From this, summation-by-parts (SBP), the discrete analogy to integration-by-parts (IBP), can be defined. 
    
    \begin{definition}[Summation-By-Parts Operator]\label{def:sbp}
        A set of operators $\mathbf{M}$, $\mathbf{D}$, $\mathbf{P}$, and $\mathbf{B}$ defined on $\mathbf{x}\in\bar{\Omega}^n$ is said to define a set of SBP operators if
        \begin{equation*}
            \mathbf{MD} + \mathbf{D}^T\mathbf{M} = \mathbf{P}^T\mathbf{BP}.
        \end{equation*} 
        The analogy to IBP can be seen through pre-multiplication by $\mathbf{v}^T$ and post-multiplication by $\mathbf{w}$.
        \begin{equation*}
            \mathbf{v}^T\mathbf{MDw} + \mathbf{v}^T\mathbf{D}^T\mathbf{Mw} = \mathbf{v}^T\mathbf{P}^T\mathbf{BPw}.
        \end{equation*}
    \end{definition}
    By setting $\mathbf{v}=\mathbf{1} = [1, \ldots , 1]^T$, the approximate integral of the derivative of the function $w$ can be calculated as 
    \begin{equation}
        \mathbf{1}^T\mathbf{MDw} = \mathbf{BPw} = \mathbf{w}_R - \mathbf{w}_L.
    \end{equation}

\section{Riemann Difference Scheme}\label{sec:scheme}
    To introduce the Riemann difference~(RD) scheme, consider the one dimensional conservation law of $m$ variables for a solution $\mathbf{u}$ and flux $\mathbf{f}$,
    \begin{equation}
        \pxvar{\mathbf{u}}{t} + \pxvar{\mathbf{f}}{x} = 0,
    \end{equation}
    where the domain $(x,t)\in\Omega\times\mathbb{R}_+$ is partitioned into $N$ elements $\Omega_k$ such that $\Omega = \bigcup_N\Omega_k$ and $\Omega_i\cap\Omega_j=\emptyset$ for $i\neq j$. These elements are then discretized with two sets of points similarly to the spectral difference method proposed by \citet{Liu2004} but first presented as the staggered grid method of \citet{Kopriva1996}. For a polynomial approximation of order $p$, the solution in each element is defined on a set of $p+1$ points $\{x_1,\dots,x_{p+1}\}\in\Omega_k$, and the flux is defined on a set of $p+2$ points $\{x_{1+\shalf},\dots,x_{p+\shalf}\}\in\Omega_k$ where $x_{i+\shalf}=\frac{1}{2} (x_i+x_{i+1})$ and $\{x_{\shalf},x_{p+\sfrac{3}{2}}\}\in\partial\Omega_k$.
    
    \begin{figure}[tbhp]
        \centering
        \adjustbox{width=0.5\linewidth}{\begin{tikzpicture}[every node/.style={font=\small},scale=2]
    
    \begin{scope}
        \draw[thick] (-1,0) -- (1,0);
        \draw[thick] (-1,-0.15) -- (-1, 0.15);
        \draw[thick] ( 1,-0.15) -- ( 1, 0.15);
        
        \draw[thick,densely dotted] (-1.1,-0.2) -- (-1.1, 1);
        \draw[thick,densely dotted] ( 1.1,-0.2) -- ( 1.1, 1);
        
        \draw[thick] (-1.2,0) -- (-1.9,0);
        \draw[thick,densely dotted] (-1.9,0) -- (-2.05,0);
        \draw[thick] ( 1.2,0) -- ( 1.9,0);
        \draw[thick,densely dotted] ( 1.9,0) -- ( 2.05,0);
        \draw[thick] (-1.2,-0.15) -- (-1.2, 0.15);
        \draw[thick] ( 1.2,-0.15) -- ( 1.2, 0.15);
    \end{scope}
    
    \begin{scope}[on above layer]
        \draw[fill=Set1-B] (-0.775,0.3) circle (0.05);
        \draw[fill=Set1-B] ( 0.000,0.2) circle (0.05);
        \draw[fill=Set1-B] ( 0.775,0.25) circle (0.05);
        
        \filldraw[Set1-A] (-1.05,0.85) -- (-1.15,0.95) -- (-1.05,0.95) -- cycle;
        \filldraw[Set1-D] (-1.05,0.85) -- (-1.15,0.95) -- (-1.15,0.85) -- cycle;
        \draw[] (-1.15,0.9-0.05) rectangle ++(0.1,0.1);
        
        \draw[fill=Set1-A] (-0.437,0.55-0.05) rectangle ++(0.1,0.1);
        \draw[fill=Set1-A] ( 0.337,0.45-0.05) rectangle ++(0.1,0.1);
        
        \filldraw[Set1-D] (1.05,0.55) -- (1.15,0.55) -- (1.15,0.65) -- cycle;
        \filldraw[Set1-A] (1.05,0.55) -- (1.15,0.65) -- (1.05,0.65) -- cycle;
        \draw[] (1.05,0.6-0.05) rectangle ++(0.1,0.1);
    \end{scope}
    
    \begin{scope}[on above layer]
        \draw[fill=Set1-C] (-2.2+0.775,0.225) circle (0.05);
        \draw[fill=Set1-D] (-1.813-0.05,0.8-0.05) rectangle ++(0.1,0.1);
        
        \draw[fill=Set1-C] (2.2-0.775,0.15) circle (0.05);
        \draw[fill=Set1-D] ( 1.813-0.05,0.4-0.05) rectangle ++(0.1,0.1);
    \end{scope}
    
    \begin{scope}[on behind layer]
        \draw[thick,densely dashed] (-0.775,0.3) -- (-0.775,0);
        \draw[thick,densely dashed] ( 0.000,0.2) -- ( 0.000,0);
        \draw[thick,densely dashed] ( 0.775,0.25) -- ( 0.775,0);
        \draw[thick] plot[smooth,domain=-1:1] (\x,{0.5*(0.24974*(\x)^2 - 0.0645161*\x + 0.4)});
        
        \draw[thick,densely dashed] (-0.775,0.3) -- (-0.775,0);
        \draw[thick,densely dashed] ( 0.000,0.2) -- ( 0.000,0);
        \draw[thick,densely dashed] ( 0.775,0.25) -- ( 0.775,0);
        \draw[thick] plot[smooth,domain=-1:1] (\x,{0.5*(0.24974*(\x)^2 - 0.0645161*\x + 0.4)});
        
        \draw[thick] plot[smooth,domain=-1:1] (\x,{-0.0244652*(\x)^3 + 0.294038*(\x)^2 - 0.125535*\x + 0.455962});
    \end{scope}
    
    \begin{scope}[on behind layer]
        \draw[thick,densely dashed] (-2.2+0.775,0) -- (-2.2+0.775,0.225);
        \draw[thick] plot[smooth,domain=0.3:1] ({\x-2.2},{0.270552*(\x)^2 - 0.0483871*\x + 0.1});
        \draw[thick,densely dotted] plot[smooth,domain=0.15:0.3] ({\x-2.2},{0.270552*(\x)^2 - 0.0483871*\x + 0.1});

        \draw[thick,densely dashed] (2.2-0.775,0) -- (2.2-0.775,0.15);
        \draw[thick] plot[smooth,domain=-1:-0.3] ({\x+2.2},{0.145682*(\x)^2 + 0.0483871*\x + 0.1});
        \draw[thick,densely dotted] plot[smooth,domain=-0.3:-0.15] ({\x+2.2},{0.145682*(\x)^2 + 0.0483871*\x + 0.1});

        \draw[thick] plot[smooth,domain=0.3:1] ({\x-2.2},{-0.122554*(\x)^3 + 0.147019*(\x)^2 + 0.147554*\x + 0.727981});
        \draw[thick,densely dotted] plot[smooth,domain=0.15:0.3] ({\x-2.2},{-0.122554*(\x)^3 + 0.147019*(\x)^2 + 0.147554*\x + 0.727981});

        \draw[thick] plot[smooth,domain=-1:-0.3] ({\x+2.2},{-0.303915*(\x)^3 + 0.117615*(\x)^2 + 0.303915*\x + 0.482385});
        \draw[thick,densely dotted] plot[smooth,domain=-0.3:-0.15] ({\x+2.2},{-0.303915*(\x)^3 + 0.117615*(\x)^2 + 0.303915*\x + 0.482385});
    \end{scope}
    
    \node[black] at (-1.425,-0.2) {$(0)$};
    
    \node[black] at (-0.775,-0.2) {$1$};
    \node[black] at ( 0.0  ,-0.2) {$2$};
    \node[black] at ( 0.775,-0.2) {$3$};
    
    \node[black] at ( 1.425,-0.2) {$(4)$};
    
    \node[black] at (-1.000,-0.4) {$\shalf$};
    \node[black] at (-0.387,-0.4) {$1+\shalf$};
    \node[black] at ( 0.387,-0.4) {$2+\shalf$};
    \node[black] at ( 1.000,-0.4) {$3+\shalf$};
    
    \begin{scope}[on behind layer]
        \draw[thick] (-0.387,0.55) circle (0.1);
        \draw[thick] (-0.487,0.55) -- ++(-0.40,0.693);
        \draw[thick] (-0.287,0.55) -- ++( 0.40,0.693);
        \draw[thick] (-0.387,1.5) circle (0.563);
        
        \draw[thick] (-0.860,1.8) -- (-0.387,1.8);
        \draw[thick] (-0.387,1.3) -- ( 0.143,1.3);
        \draw[thick, dotted] (-0.387,1.8) -- (-0.387,1.3);
        
        \draw[fill=Set1-B] (-0.775,1.8) circle (0.05);
        \draw[fill=Set1-B] ( 0.000,01.3) circle (0.05);
        \draw[fill=Set1-A] (-0.387,1.5) circle (0.05);
    \end{scope}
    
\end{tikzpicture}}
        \caption{\label{fig:schema_1d}Diagram of the 1D RD scheme for $\mathbb{P}_2$. Circles represent the solution points and squares represent the flux points.}
    \end{figure}
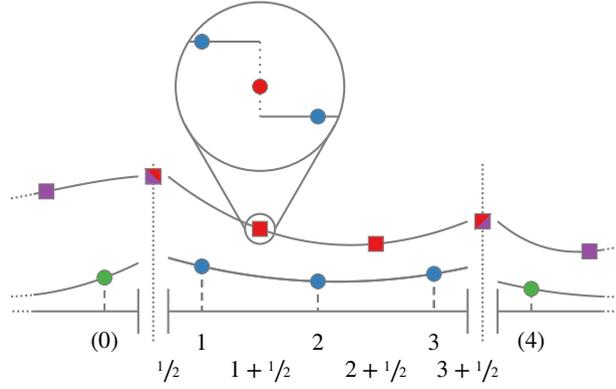
    
    Through these two sets of points, polynomial interpolations for the solution and the flux in the sub-domain $\Omega_k$ can be formed as
    \begin{equation}
        \hat{\mathbf{u}}_k = \sum^{p+1}_{j=1} \mathbf{u}(x_j)l^u_{j}(x) \quad \mathrm{and} \quad \hat{\mathbf{f}}_k = \sum^{p+1}_{j=0}\mathbf{f}\big(x_{j+\shalf}\big)l^f_{j+\shalf}(x),
    \end{equation}
    where $l^u$ and $l^f$ are the Lagrange nodal basis functions defined as
    \begin{equation}
        l^u_i(x) = \prod^{p+1}_{\substack{j=1\\j\neq i}}\frac{x-x_j}{x_i-x_j} \quad \mathrm{and} \quad l^f_{i+\shalf}(x) = \prod^{p+1}_{\substack{j=0\\j\neq i}}\frac{x-x_{j+\shalf}}{x_{i+\shalf}-x_{j+\shalf}}.
    \end{equation}
    For brevity, we drop the subscript $k$ and present the scheme for an arbitrary $\Omega_k$. We then define the term
    \begin{equation}
        c_{ij+\shalf} = \dx{l^f_{j+\shalf}(x)}{x}\bigg|_{x=x_i},
    \end{equation}
    from which it can be seen that 
    \begin{equation}\label{eq:dfdx}
        \dx{\hat{\mathbf{f}}}{x}\bigg|_{x=x_i} = \sum^{p+1}_{j=0}c_{ij+\shalf}\mathbf{f}\big(x_{j+\shalf}\big).
    \end{equation}
    It can be easily shown that $c_{ij+\shalf}$ has the property that $\sum^{p+1}_{j=0} c_{ij+\shalf}=0$.
    
    In contrast to the spectral difference (SD) scheme, the Riemann difference scheme utilizes an auxiliary (Lax--Friedrichs) flux at the flux points and, as such, takes the semi-discrete form for the $i$-th solution point as
    \begin{equation}\label{eq:rd}
        \pxvar{\mathbf{u}_i}{t} = -\sum^{p+1}_{j=0}c_{ij+\shalf}\ob{f}_{ij+\shalf}
    \end{equation}
    where the auxiliary flux is set as
    \begin{equation}
        \ob{f}_{ij+\shalf} = \half(\mathbf{f}_{j}+\mathbf{f}_{j+1}) - \frac{d_{ij+\shalf}}{2c_{ij+\shalf}}(\mathbf{u}_{j+1} - \mathbf{u}_{j})
    \end{equation}
    for some scalar $d_{ij+\shalf}>0$. Additionally, we define the auxiliary state as 
    \begin{equation}
        \ob{u}_{ij+\shalf} = \half(\mathbf{u}_{j}+\mathbf{u}_{j+1}) - \frac{c_{ij+\shalf}}{2d_{ij+\shalf}}(\mathbf{f}_{j+1} - \mathbf{f}_{j}).
    \end{equation}
    By setting $d_{ij+\shalf}$ as
    \begin{equation}\label{eq:dij}
        d_{ij+\shalf} = \lambda_{\max} (\mathbf{u}_j, \mathbf{u}_{j+1}) |c_{ij+\shalf} |,
    \end{equation}
    where $\lambda_{\max} (\mathbf{u}_j, \mathbf{u}_{j+1})$ denotes the maximum wave speed defined in \cref{ssec:invariant_sets} for a Riemann problem with $\mathbf{u}_l = \mathbf{u}_j, \mathbf{u}_r = \mathbf{u}_{j+1}$, it can be seen that $\ob{f}_{ij+\shalf}$ forms a Lax--Friedrichs flux and $\ob{u}_{ij+\shalf}$ takes the form of $\ob{v}$ in \cref{rieman_lemma}.
    
    The auxiliary state and auxiliary flux may be related by
    \begin{equation}\label{eq:fbar_ubar}
       c_{ij+\shalf}\ob{f}_{ij+\shalf}=-d_{ij+\shalf}\ob{u}_{ij+\shalf} + d_{ij+\shalf}\mathbf{u}_{j} + c_{ij+\shalf}\mathbf{f}_j.
    \end{equation}
    In both the auxiliary state and auxiliary flux, we use the convention that $u_0$ is the nearest solution point in the element to the left and $u_{p+2}$ is the nearest solution point in the element to the right. A diagram of the scheme and point layout is shown in \cref{fig:schema_1d}.
    
    If we then apply forward Euler temporal integration, the temporal update can be written as
    \begin{equation}\label{eq:rd_euler}
        \mathbf{u}^{n+1}_i = \mathbf{u}^{n}_i - \Delta t\sum^{p+1}_{j=0}c_{ij+\shalf}\ob{f}_{ij+\shalf}^{n} = \mathbf{u}^{n}_i - \Delta t\sum^{p+1}_{j=0} \left ( -d_{ij+\shalf}\ob{u}_{ij+\shalf}^{n} + d_{ij+\shalf}\mathbf{u}_{j}^{n} + c_{ij+\shalf}\mathbf{f}_j^{n} \right ).
    \end{equation}
    where $t^{n+1} = t^{n} + \Delta t, n \geqslant 0$. Here we utilize forward Euler due to its simplicity and strong stability preservation~(SSP) property, but the theoretical results of this section may be extended to SSP Runge--Kutta~(SSP-RK) schemes. From \cref{eq:rd_euler}, it is straightforward to extend the scheme over $\Omega$.
    
    We may now move on to state and prove the properties of the scheme. We will subsequently utilize the convention that $(\cdot)_{i,k}^n$ denotes the value of the $i$-th point of $(\cdot)$ in $\Omega_k$ at time step $n$.
    
    \begin{theorem}[Conservation]\label{thm:conservation} 
        Let $\mathbf{U}_k^n$ be the solution defined by the solution points $\{\mathbf{u}^n_{1,k},\dots,\mathbf{u}^n_{p+1,k}\}$. The scheme defined by \cref{eq:rd_euler} is conservative in the sense that for all $n \geqslant 0$,
        \begin{equation*}
            \sum_{k=1}^{N} \langle\mathbf{U}^{n+1}_{k},\mathbf{1}\rangle = \sum_{k=1}^{N} \langle\mathbf{U}^{n}_{k},\mathbf{1}\rangle - \Delta t(\ob{f}_{R,N} - \ob{f}_{L,1}),
        \end{equation*}
        where $\langle \cdot, \cdot \rangle$ denotes the inner product and $\ob{f}_{R,N}$ and $\ob{f}_{L,1}$ are the incoming and outgoing Riemann fluxes, respectively, at the domain boundaries.
    \end{theorem}
    
    \newproof{pot_cons}{Proof of \cref{thm:conservation}}
    \begin{pot_cons}
        By introducing the quadrature
        \begin{equation}
            \int_{\Omega_k}\mathbf{U}^n_{k}\ \mathrm{d}\mathbf{x} \approx \langle\mathbf{U}^{n}_{k},\mathbf{1}\rangle = \sum^{p+1}_{i=1}m_i\mathbf{u}^n_{i,k}
        \end{equation}
        and substituting it into \cref{eq:rd_euler}, we obtain
        \begin{equation}
            \sum_{k=1}^{N} \sum^{p+1}_{i=1}m_i\mathbf{u}^{n+1}_{i,k} = \sum_{k=1}^{N} \sum^{p+1}_{i=1}m_i\mathbf{u}^{n}_{i,k} - \Delta t\sum_{k=1}^{N} \sum_{i=1}^{p+1}\sum^{p+1}_{j=0}m_ic_{ij+\shalf}\ob{f}^n_{ij+\shalf,k}.
        \end{equation}
        With the summation-by-parts framework and \cref{def:sbp}, it can be seen that for any arbitrary $\Omega_k$,
        \begin{equation*}
            \sum_{i=1}^{p+1}\sum^{p+1}_{j=0}m_ic_{ij+\shalf}\ob{f}_{ij+\shalf,k} = \mathbf{1}^T\mathbf{MD}\ob{f}_{:, k} = \mathbf{BP}\ob{f} = \ob{f}_{R,k} - \ob{f}_{L,k},
        \end{equation*}
        where $\ob{f}_{:, k}$ denotes the vector of Lax--Friedrichs fluxes and  $\ob{f}_{L,k}$ and $\ob{f}_{R,k}$ denote the left and right interface fluxes, respectively. As the interface fluxes take a common value for adjacent elements, i.e., $\ob{f}_{L,k} = \ob{f}_{R,k-1}$ and vice versa, the summation over all elements yields
        \begin{equation}
            \sum_{k=1}^{N} \sum_{i=1}^{p+1}\sum^{p+1}_{j=0} m_i c_{ij+\shalf}\ob{f}_{ij+\shalf,k} =  \ob{f}_{R,N} - \ob{f}_{L,1}.
        \end{equation}
        The conclusion follows from the definition of the inner product quadrature. 
        \qed
    \end{pot_cons}
    
    \begin{theorem}[Convergence]\label{thm:convergence}
        The scheme defined by \cref{eq:rd_euler} converges in the sense that
        \begin{equation*}
            \lim_{\substack{\Delta t, \Delta x \rightarrow 0}}\|\mathbf{u}_i - \mathbf{u}(x_i)\|_2 = 0.
        \end{equation*}
    \end{theorem}
    
    \newproof{pot_conv}{Proof of \cref{thm:convergence}}
    \begin{pot_conv}
        Let $C_t,C_h\in\mathbb{R}_+$. The forwards Euler approximation of the temporal derivative gives
        \begin{equation*}
            \mathbf{u}_i^{n+1} = \mathbf{u}_i^n + \Delta t\pxvar{\ob{f}}{x}^n + C_t\Delta t^{2}.
        \end{equation*}
        From \citet{Lax1957}, it is known that $\ob{f}$ provides the flux with accuracy $\mathcal{O}(\Delta x^2)$. Let $\mathbf{e}(x_{i+\shalf})$ be the error in the flux at the flux points $x_{i+\shalf}$. This produces a polynomial, $\hat{\mathbf{e}}(x)$, of degree $p+1$, which is exactly differentiated by $\mathbf{D}$. However, $\sup_i \hat{\mathbf{e}}(\mathbf{x}_i)=\mathcal{O}(\Delta x^2)$, which implies $\lim_{\Delta x\rightarrow 0}\|\pxvar{\hat{\mathbf{e}}(x)}{x}\|_\infty = \lim_{\Delta x\rightarrow 0}C_h\Delta x = 0$. As a result, 
        \begin{equation*}
            \lim_{\substack{\Delta t, \Delta x \rightarrow 0}}\|\mathbf{u}_i - \mathbf{u}(x_i)\|_2 =  \lim_{\substack{\Delta t, \Delta x \rightarrow 0}}\left(C_t\Delta t^{2} + C_h\Delta x^{1}\right) = 0. \quad \hfill \quad \qed
        \end{equation*}   
    \end{pot_conv}
    
    \begin{theorem}[Local Invariance]\label{thm:local_invariance} 
         For some sub-domain $\Omega_k$, let $\mathbf{u}_{i}^n = \mathbf{u}_{i,k}^n \in\mathcal{B}\ \forall\ i \in \{1,...,p+1\}$ for the set $\mathcal{B}$ defined in \cref{def:invariant_set}. For the scheme defined by \cref{eq:rd_euler} and $d_{ij+\shalf}$ defined by \cref{eq:dij}, there exists a strictly positive $\Delta t$ such that $\mathbf{u}^{n+1}_i \in \mathcal{B}$ if $G_\mathcal{B}(\mathbf{u}^{n}_i), G_\mathcal{B}(\ob{u}^{n}_{ii+\shalf}) \neq 0$.
        
    \end{theorem}
    
    \newproof{pot_li}{Proof of \cref{thm:local_invariance}}
    
    \begin{pot_li}
        By substituting \cref{eq:fbar_ubar} into \cref{eq:rd_euler}, we obtain
        \begin{equation}
            \mathbf{u}^{n+1}_i = \mathbf{u}^{n}_i - \Delta t \sum^{p+1}_{j=0}\left(d_{ij+\shalf}\mathbf{u}^{n}_j + c_{ij+\shalf}\mathbf{f}^{n}_j\right) + \Delta t\sum_{j=0}^{p+1} d_{ij+\shalf}\ob{u}_{ij+\shalf}
        \end{equation}
        which may be rewritten as         
        \begin{equation}\label{eq:comb_form}
            \mathbf{u}^{n+1}_i = \left[1 - \Delta t d_{ii+\shalf} \right] \mathbf{u}^{n}_i - \Delta t \sum^{p+1}_{\substack{j=0 \\ j\neq i}}\left(d_{ij+\shalf}\mathbf{u}^{n}_j + c_{ij+\shalf}\mathbf{f}^{n}_j + \frac{c_{ii+\shalf}}{p+1} \mathbf{f}^{n}_i\right) + \Delta t\sum_{j=0}^{p+1} d_{ij+\shalf}\ob{u}^{n}_{ij+\shalf}.
        \end{equation}
        We then define $\mathbf{v}_{ij}$ as
        \begin{equation}
            \mathbf{v}^{n}_{ij} = \ob{u}^{n}_{ij+\shalf} - \left(\mathbf{u}^{n}_j +  \frac{c_{ij+\shalf}}{d_{ij+\shalf}}\mathbf{f}^{n}_j +  \frac{c_{ii+\shalf}}{(p+1)d_{ii+\shalf}} \mathbf{f}^{n}_i\right) =   \frac{c_{ii+\shalf}}{(p+1)d_{ii+\shalf}} \mathbf{f}_i^n - \frac{c_{ij+\shalf}}{d_{ij+\shalf}} \ob{f}^n_{ij+\shalf}.
        \end{equation}
        It can be seen that this term is not necessarily in $\mathcal{B}$ but allows \cref{eq:comb_form} to be expressed as 
        \begin{equation}
            \mathbf{u}^{n+1}_i = \left[1 - \Delta t d_{ii+\shalf}\right] \mathbf{u}^{n}_i + \Delta t \sum^{p+1}_{\substack{j=0 \\ j\neq i}} d_{ij+\shalf} \mathbf{v}^{n}_{ij} + \Delta t d_{ii+\shalf}\ob{u}^{n}_{ii+\shalf} = a_1 \mathbf{w}_1 + a_2 \mathbf{w}_2 + a_3 \mathbf{w}_3
        \end{equation}
        where $a_1 = \left[1 - \Delta t d_{ii+\shalf}\right]$, $a_2 = \Delta t$, and $a_3 = \Delta t d_{ii+\shalf}$. Furthermore, the states $\mathbf{w}_1, \mathbf{w}_3$ are in $\mathcal{B}$ while $\mathbf{w}_2$ is not necessarily in $\mathcal{B}$ \citep{Guermond2019}.  As $d_{ii+\shalf}$ is strictly positive, there always exists some strictly positive $\Delta t \leqslant \frac{1}{d_{ii+\shalf}}$ such that the state $\mathbf{w}^* = a_1 \mathbf{w}_1 + a_3 \mathbf{w}_3$ forms a convex combination of states in $\mathcal{B}$. As a result, $\mathbf{w}^* \in \mathcal{B}$.
        
        From \cref{cor:noninvariant_sum}, $\mathbf{u}^{n+1}_i =  \mathbf{w}^* + a_2 \mathbf{w_2} \in \mathcal{B}$ if
        \begin{equation*}
            a_2 = \Delta t \leqslant \frac{G_\mathcal{B}(\mathbf{w}^*)}{\|\mathbf{w}_2\|_\mathcal{B}}.
        \end{equation*}
        Therefore, from \cref{lem:convex_dist}, there exists a $\Delta t$ such that $\mathbf{u}^{n+1}_i \in \mathcal{B}$ if
        \begin{equation*}
            \Delta t \leqslant \min \left [ \frac{1}{d_{ii+\shalf}},  \frac{G_\mathcal{B}(\mathbf{u}^{n}_i)}{\|\mathbf{w}_2\|_\mathcal{B}}, \frac{G_\mathcal{B}(\ob{u}^{n}_{ii+\shalf})}{\|\mathbf{w}_2\|_\mathcal{B}} \right].
        \end{equation*}
        By extension, this implies that there exists a strictly positive $\Delta t$ such that $\mathbf{u}^{n+1}_i \in \mathcal{B}$ if $G_\mathcal{B}(\mathbf{u}^{n}_i), G_\mathcal{B}(\ob{u}^{n}_{ii+\shalf}) \neq 0$.
        \qed
    \end{pot_li}
    
    \begin{corollary}[Global Invariance]\label{cor:global_invariance} 
        Let $I = \{1,...,p+1\}$, $K =  \{1,...,N\}$, and $\mathcal{B}$ be the set defined in \cref{def:invariant_set}. If $\mathbf{u}_{i,k}^n \in\mathcal{B}$ and $G_\mathcal{B}(\mathbf{u}^{n}_{i,k}), G_\mathcal{B}(\ob{u}^{n}_{ii+\shalf,k}) \neq 0$ for all $i \in I$ and $k \in K$, then there exists a strictly positive $\Delta t$ such that $\mathbf{u}^{n+1}_{i,k} \in \mathcal{B}$ for all $i \in I$ and $k \in K$, i.e., the RD scheme provides a mapping such that $\mathcal{B}$ is an invariant domain.
    \end{corollary}

\section{Numerical Implementation}\label{sec:numerical}
    Although it is shown that the solution produced by the RD method is invariant domain preserving under set conditions, these benefits in numerical stability comes at the expense of order of accuracy. It is known that for hyperbolic conservation laws, the solution, in a weak sense, is smooth except for a countable number of discontinuities~\citep{Hopf1950}. Therefore, to increase the utility of the method in scale-resolving simulations, it is beneficial to couple it to a higher-order collocation method. 
    The method of choice in this work is the flux reconstruction (FR) scheme of \citet{Huynh2007,Vincent2010}, implemented within the PyFR software package \citep{Witherden2014}. To control the switching between the schemes, the sensor of \citet{Persson2006} was used; however, alternative choices of sensors could be more appropriate in cases where the use of tunable parameters is not feasible. In the subsequent sections, the use of the FR scheme paired with the RD scheme via a sensor is denoted by RD-FR whereas the use of the schemes independently of each other is denoted by FR and RD, respectively. 
    
    \begin{figure}[tbhp]
        \centering
        \adjustbox{width=0.35\linewidth}{\begin{tikzpicture}[every node/.style={font=\small},scale=2]

    \begin{scope}
        \draw[thick](-1,-1) rectangle ++(2,2);
        
        \draw[thick] (-1.9,-1) -- (-1.2,-1) -- (-1.2, 1) -- (-1.9, 1);
        \draw[thick,densely dotted] (-1.9,-1) -- (-2.05,-1);
        \draw[thick,densely dotted] (-1.9, 1) -- (-2.05, 1);
        
        \draw[thick] ( 1.9,-1) -- ( 1.2,-1) -- ( 1.2, 1) -- ( 1.9, 1);
        \draw[thick,densely dotted] ( 1.9,-1) -- ( 2.05,-1);
        \draw[thick,densely dotted] ( 1.9, 1) -- ( 2.05, 1);
        
        \draw[thick] (-1, 1.9) -- (-1, 1.2) -- (1, 1.2) -- (1, 1.9);
        \draw[thick,densely dotted] (-1, 1.9) -- (-1, 2.05);
        \draw[thick,densely dotted] ( 1, 1.9) -- ( 1, 2.05);
        
        \draw[thick] (-1,-1.9) -- (-1,-1.2) -- (1,-1.2) -- (1,-1.9);
        \draw[thick,densely dotted] (-1,-1.9) -- (-1,-2.05);
        \draw[thick,densely dotted] ( 1,-1.9) -- ( 1,-2.05);
        
        \draw[thick,densely dotted] (-1, 1.1) -- ( 1, 1.1);
        \draw[thick,densely dotted] (-1,-1.1) -- ( 1,-1.1);
        \draw[thick,densely dotted] ( 1.1,-1) -- ( 1.1, 1);
        \draw[thick,densely dotted] (-1.1,-1) -- (-1.1, 1);
    \end{scope}
    \begin{scope}[on above layer]
        \draw[fill=Set1-B, opacity=1.0] (-0.775,-0.775) circle (0.05);
        \draw[fill=Set1-B, opacity=1.0] ( 0.000,-0.775) circle (0.05);
        \draw[fill=Set1-B] ( 0.775,-0.775) circle (0.05);
        \draw[fill=Set1-B] (-0.775, 0.000) circle (0.05);
        \draw[fill=Set1-B] ( 0.000, 0.000) circle (0.05);
        \draw[fill=Set1-B] ( 0.775, 0.000) circle (0.05);
        \draw[thick,rounded corners=4pt] (0.775-0.12,0.000-0.12) rectangle ++(0.24,0.24) node[below,xshift=-0.8em,yshift=-1.3em,black] {\normalsize $i$};
        
        \draw[fill=Set1-B, opacity=1.0] (-0.775, 0.775) circle (0.05);
        \draw[fill=Set1-B, opacity=1.0] ( 0.000, 0.775) circle (0.05);
        \draw[fill=Set1-B] ( 0.775, 0.775) circle (0.05);
    \end{scope}
       
    \begin{scope}[on above layer]
        \draw[fill=Set1-A, opacity=1.0] (-0.437,-0.825) rectangle ++(0.1,0.1);
        \draw[fill=Set1-A, opacity=1.0] ( 0.337,-0.825) rectangle ++(0.1,0.1);
        \draw[fill=Set1-A] (-0.437,-0.050) rectangle ++(0.1,0.1);
        \draw[fill=Set1-A] ( 0.337,-0.050) rectangle ++(0.1,0.1);
        \draw[fill=Set1-A, opacity=1.0] (-0.437, 0.725) rectangle ++(0.1,0.1);
        \draw[fill=Set1-A, opacity=1.0] ( 0.337, 0.725) rectangle ++(0.1,0.1);
        
        \draw[fill=Set1-A, opacity=1.0] (-0.825,-0.437) rectangle ++(0.1,0.1);
        \draw[fill=Set1-A, opacity=1.0] (-0.825, 0.337) rectangle ++(0.1,0.1);
        \draw[fill=Set1-A, opacity=1.0] (-0.050,-0.437) rectangle ++(0.1,0.1);
        \draw[fill=Set1-A, opacity=1.0] (-0.050, 0.337) rectangle ++(0.1,0.1);
        \draw[fill=Set1-A] ( 0.725,-0.437) rectangle ++(0.1,0.1);
        \draw[fill=Set1-A] ( 0.725, 0.337) rectangle ++(0.1,0.1);
   \end{scope}
   
   \begin{scope}[on above layer]
        \filldraw[Set1-A, opacity=1.0] (-1.05,0.775-0.05) -- (-1.15,0.775+0.05) -- (-1.05,0.775+0.05) -- cycle;
        \filldraw[Set1-D, opacity=1.0] (-1.05,0.775-0.05) -- (-1.15,0.775+0.05) -- (-1.15,0.775-0.05) -- cycle;
        \draw[] (-1.15,0.775-0.05) rectangle ++(0.1,0.1);
        
        \filldraw[Set1-A] (-1.05,0.000-0.05) -- (-1.15,0.000+0.05) -- (-1.05,0.000+0.05) -- cycle;
        \filldraw[Set1-D] (-1.05,0.000-0.05) -- (-1.15,0.000+0.05) -- (-1.15,0.000-0.05) -- cycle;
        \draw[] (-1.15,0.000-0.05) rectangle ++(0.1,0.1);
        
        \filldraw[Set1-A, opacity=1.0] (-1.05,-0.775-0.05) -- (-1.15,-0.775+0.05) -- (-1.05,-0.775+0.05) -- cycle;
        \filldraw[Set1-D, opacity=1.0] (-1.05,-0.775-0.05) -- (-1.15,-0.775+0.05) -- (-1.15,-0.775-0.05) -- cycle;
        \draw[] (-1.15,-0.775-0.05) rectangle ++(0.1,0.1);
        
        \filldraw[Set1-A, opacity=1.0] ( 1.05, 0.775-0.05) -- ( 1.15, 0.775+0.05) -- ( 1.05, 0.775+0.05) -- cycle;
        \filldraw[Set1-D, opacity=1.0] ( 1.05, 0.775-0.05) -- ( 1.15, 0.775+0.05) -- ( 1.15, 0.775-0.05) -- cycle;
        \draw[] ( 1.05, 0.775-0.05) rectangle ++(0.1,0.1);
        
        \filldraw[Set1-A] ( 1.05, 0.000-0.05) -- ( 1.15, 0.000+0.05) -- ( 1.05, 0.000+0.05) -- cycle;
        \filldraw[Set1-D] ( 1.05, 0.000-0.05) -- ( 1.15, 0.000+0.05) -- ( 1.15, 0.000-0.05) -- cycle;
        \draw[] ( 1.05, 0.000-0.05) rectangle ++(0.1,0.1);
        
        \filldraw[Set1-A, opacity=1.0] ( 1.05,-0.775-0.05) -- ( 1.15,-0.775+0.05) -- ( 1.05,-0.775+0.05) -- cycle;
        \filldraw[Set1-D, opacity=1.0] ( 1.05,-0.775-0.05) -- ( 1.15,-0.775+0.05) -- ( 1.15,-0.775-0.05) -- cycle;
        \draw[] ( 1.05,-0.775-0.05) rectangle ++(0.1,0.1);
        
        \filldraw[Set1-A] ( 0.775-0.05, 1.05) -- ( 0.775+0.05, 1.15) -- ( 0.775+0.05, 1.05) -- cycle;
        \filldraw[Set1-D] ( 0.775-0.05, 1.05) -- ( 0.775+0.05, 1.15) -- ( 0.775-0.05, 1.15) -- cycle;
        \draw[] ( 0.775-0.05, 1.05) rectangle ++(0.1,0.1);
        
        \filldraw[Set1-A, opacity=1.0] ( 0.000-0.05, 1.05) -- ( 0.000+0.05, 1.15) -- ( 0.000+0.05, 1.05) -- cycle;
        \filldraw[Set1-D, opacity=1.0] ( 0.000-0.05, 1.05) -- ( 0.000+0.05, 1.15) -- ( 0.000-0.05, 1.15) -- cycle;
        \draw[] ( 0.000-0.05, 1.05) rectangle ++(0.1,0.1);
        
        \filldraw[Set1-A, opacity=1.0] (-0.775-0.05, 1.05) -- (-0.775+0.05, 1.15) -- (-0.775+0.05, 1.05) -- cycle;
        \filldraw[Set1-D, opacity=1.0] (-0.775-0.05, 1.05) -- (-0.775+0.05, 1.15) -- (-0.775-0.05, 1.15) -- cycle;
        \draw[] (-0.775-0.05, 1.05) rectangle ++(0.1,0.1);
        
        \filldraw[Set1-A] ( 0.775-0.05,-1.05) -- ( 0.775+0.05,-1.15) -- ( 0.775+0.05,-1.05) -- cycle;
        \filldraw[Set1-D] ( 0.775-0.05,-1.05) -- ( 0.775+0.05,-1.15) -- ( 0.775-0.05,-1.15) -- cycle;
        \draw[] ( 0.775-0.05,-1.15) rectangle ++(0.1,0.1);
        
        \filldraw[Set1-A, opacity=1.0] ( 0.000-0.05,-1.05) -- ( 0.000+0.05,-1.15) -- ( 0.000+0.05,-1.05) -- cycle;
        \filldraw[Set1-D, opacity=1.0] ( 0.000-0.05,-1.05) -- ( 0.000+0.05,-1.15) -- ( 0.000-0.05,-1.15) -- cycle;
        \draw[] ( 0.000-0.05,-1.15) rectangle ++(0.1,0.1);
        
        \filldraw[Set1-A, opacity=1.0] (-0.775-0.05,-1.05) -- (-0.775+0.05,-1.15) -- (-0.775+0.05,-1.05) -- cycle;
        \filldraw[Set1-D, opacity=1.0] (-0.775-0.05,-1.05) -- (-0.775+0.05,-1.15) -- (-0.775-0.05,-1.15) -- cycle;
        \draw[] (-0.775-0.05,-1.15) rectangle ++(0.1,0.1);
    \end{scope}

    \begin{scope}[on above layer]
        \draw[fill=Set1-C, opacity=1.0] (-2.2+0.775,-0.775) circle (0.05);
        \draw[fill=Set1-C] (-2.2+0.775, 0.000) circle (0.05);
        \draw[fill=Set1-C, opacity=1.0] (-2.2+0.775, 0.775) circle (0.05);
        
        \draw[fill=Set1-D, opacity=1.0] (-1.813-0.05, 0.775-0.05) rectangle ++(0.1,0.1);
        \draw[fill=Set1-D] (-1.813-0.05, 0.000-0.05) rectangle ++(0.1,0.1);
        \draw[fill=Set1-D, opacity=1.0] (-1.813-0.05,-0.775-0.05) rectangle ++(0.1,0.1);
        
        \draw[fill=Set1-D, opacity=1.0] (-1.425-0.05,-0.387-0.05) rectangle ++(0.1,0.1);
        \draw[fill=Set1-D, opacity=1.0] (-1.425-0.05, 0.387-0.05) rectangle ++(0.1,0.1);
        
        \draw[fill=Set1-C, opacity=1.0] ( 2.2-0.775,-0.775) circle (0.05);
        \draw[fill=Set1-C] ( 2.2-0.775, 0.000) circle (0.05);
        \draw[fill=Set1-C, opacity=1.0] ( 2.2-0.775, 0.775) circle (0.05);
        
        \draw[fill=Set1-D, opacity=1.0] ( 1.813-0.05, 0.775-0.05) rectangle ++(0.1,0.1);
        \draw[fill=Set1-D] ( 1.813-0.05, 0.000-0.05) rectangle ++(0.1,0.1);
        \draw[fill=Set1-D, opacity=1.0] ( 1.813-0.05,-0.775-0.05) rectangle ++(0.1,0.1);
        
        \draw[fill=Set1-D, opacity=1.0] ( 1.425-0.05,-0.387-0.05) rectangle ++(0.1,0.1);
        \draw[fill=Set1-D, opacity=1.0] ( 1.425-0.05, 0.387-0.05) rectangle ++(0.1,0.1);
        
        \draw[fill=Set1-C, opacity=1.0] (-0.775,-2.2+0.775) circle (0.05);
        \draw[fill=Set1-C, opacity=1.0] ( 0.000,-2.2+0.775) circle (0.05);
        \draw[fill=Set1-C] ( 0.775,-2.2+0.775) circle (0.05);
        
        \draw[fill=Set1-D] ( 0.775-0.05,-1.813-0.05) rectangle ++(0.1,0.1);
        \draw[fill=Set1-D, opacity=1.0] ( 0.000-0.05,-1.813-0.05) rectangle ++(0.1,0.1);
        \draw[fill=Set1-D, opacity=1.0] (-0.775-0.05,-1.813-0.05) rectangle ++(0.1,0.1);
        
        \draw[fill=Set1-D, opacity=1.0] (-0.387-0.05,-1.425-0.05) rectangle ++(0.1,0.1);
        \draw[fill=Set1-D, opacity=1.0] ( 0.387-0.05,-1.425-0.05) rectangle ++(0.1,0.1);
        
        \draw[fill=Set1-C, opacity=1.0] (-0.775, 2.2-0.775) circle (0.05);
        \draw[fill=Set1-C, opacity=1.0] ( 0.000, 2.2-0.775) circle (0.05);
        \draw[fill=Set1-C] ( 0.775, 2.2-0.775) circle (0.05);
        
        \draw[fill=Set1-D] ( 0.775-0.05, 1.813-0.05) rectangle ++(0.1,0.1);
        \draw[fill=Set1-D, opacity=1.0] ( 0.000-0.05, 1.813-0.05) rectangle ++(0.1,0.1);
        \draw[fill=Set1-D, opacity=1.0] (-0.775-0.05, 1.813-0.05) rectangle ++(0.1,0.1);
        
        \draw[fill=Set1-D, opacity=1.0] (-0.387-0.05, 1.425-0.05) rectangle ++(0.1,0.1);
        \draw[fill=Set1-D, opacity=1.0] ( 0.387-0.05, 1.425-0.05) rectangle ++(0.1,0.1);
    \end{scope}
    
    \begin{scope}[on behind layer]
        \draw[Pastel1-E,fill=Pastel1-E,rounded corners=4pt] (0.775-0.12,-1.175) rectangle ++(0.24,2.35);
        \draw[Pastel1-E,fill=Pastel1-E,rounded corners=4pt] (-1.175,0.000-0.12) rectangle ++(2.35,0.24);
    \end{scope}

\end{tikzpicture}}
        \caption{\label{fig:schema_2d}Diagram of the 2D RD scheme for $\mathbb{P}_2$. Circles represent the solution points and squares represent the flux points.}
    \end{figure}
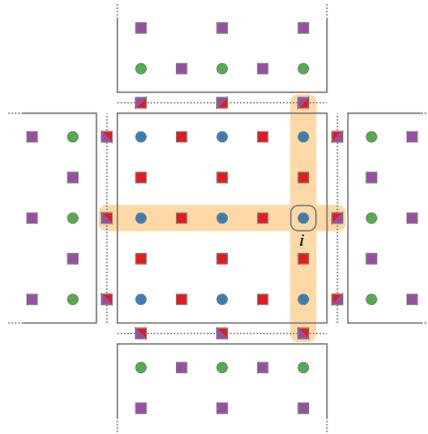
    
    The RD method is extended to higher dimensions through the use of a tensor product formulation where the gradient is calculated along lines as shown in \cref{fig:schema_2d}. This is evidently stable for affine transformations from the reference space $[-1,1]^d$, and a topic of future work is the adaptation to non-affine transformations. In the cases to be shown, some non-affine elements were used successfully.
    
    The proofs of the \cref{ssec:invariant_sets} were dependent on the strong stability of the temporal integration, e.g., forward Euler. This is also a property of the set of strong stability preserving explicit Runge--Kutta (SSP-RK) schemes of \citet{Gottlieb2001}. Hence, the properties presented in \cref{ssec:invariant_sets} follow for SSP-RK schemes, and in the following numerical cases, we utilize a three-stage, third-order SSP-RK3 scheme. The remainder of this section outlines the FR method, discontinuity sensor, and the maximum wavespeed calculation.
    
\subsection{Flux Reconstruction}\label{ssec:fr}
    The FR method may be considered as a generalization of the nodal discontinuous Galerkin method \citep{Hesthaven2008a,Zwanenburg2016}. We will give a brief description of the FR algorithm in one dimension to first-order systems, but \citet{Witherden2016} and the references therein provide details on extensions to higher dimensions and second-order PDEs. For this procedure, we utilize the reference domain $\hat{x}\in\hat{\Omega}=[-1,1]$ and  the transformation $T_k:\hat{\Omega}\rightarrow\Omega_k$ for the sub-domain $\Omega_k$. We define the reference shape functions $\hat{l_i}(\hat{x})$ as the Lagrange interpolating polynomials for a set of $p+1$ unique nodes $\{\hat{x}_0, ... , \hat{x}_p\}$, and we also define the interpolation operators $I_l$ and $I_r$ such that $I_lv = v_l = v(-1)$ and $I_rv = v_r = v(1)$ where $v\in \mathbb{P}_p$. Lastly, we define a continuous flux function $f(u)$. Therefore, for the following equation,
    \begin{equation*}
        \px{u}{t} + \px{f(u)}{x} = 0,
    \end{equation*}
    the approximation of the solution $u$ at $x_j\in\Omega_k$ is given by
    \begin{equation}
        u(x_j) = u_j \approx  \sum^p_{i=0} u(\hat{x}_i) l_i(T_k^{-1} x_j),
    \end{equation}
    and the FR spatial derivative of the flux $f$ at $x_j\in\Omega_k$ is given by
    \begin{equation}
        \px{f(x_j)}{x} \approx \left(\dx{T_k}{\hat{x}}\bigg|_{\hat{x}=T_k^{-1}x_j}\right)^{-1}\left[\sum^p_{i=0}f(u_i)\dx{l_i}{\hat{x}} + (f_l^{I} - I_l f)\dx{h_l}{\hat{x}} + (f_r^{I} - I_r f)\dx{h_r}{\hat{x}}\right]_{\hat{x}=T_k^{-1}x_j}
    \end{equation}
    Here, we have introduced the two key components of the FR method. The first is the common interface flux, denoted by $f_l^I$ and $f_r^I$ for the left and right interfaces, respectively, of the element. This is typically calculated by treating the interface as a Riemann problem using the interpolated solution at either side of interface as the initial condition. Commonly used approaches for this are approximate Riemann solvers such as that of \citet{Rusanov1962,Roe1981}. We have also introduced the correction functions $h_l, h_r \in \mathbb{P}_{p+1}$. These functions, defined in the reference domain, have the properties that $h_l(-1)=h_r(1)=1$ and $h_l(1)=h_r(-1)=0$. From this, it can be surmised that the FR algorithm 
    approximates the derivative of the flux with Lagrange interpolating polynomials and then applies a correction such that the flux approximation takes a common value at the interface which takes into account the contribution of the neighboring elements.

\subsection{Discontinuity Sensor}\label{ssec:sensor}
    In order to detect the presence of a discontinuity, the method of \citet{Persson2006} was used where the relative energy of the modal components of the solution was compared. For a two-dimensional tensor-product element, the solution $\hat{u}$ can be expressed in its modal form through the equivalence
    \begin{equation}
         \hat{u} = \sum^{p}_{i=0}\sum^{p}_{j=0} \mathbf{m}_{i,j}\mathcal{L}_{i,j},
    \end{equation}
    where $\mathcal{L}_{i,j}$ denotes the product of the $i$th and $j$th Legendre polynomials. A truncated solution, $\tilde{u}$, can then be defined as  
    \begin{equation}
         \tilde{u} = \sum^{p-1}_{i=0}\sum^{p-1}_{j=0} \mathbf{m}_{i,j}\mathcal{L}_{i,j}.
    \end{equation}
    From this, a smoothness indicator $S_e$ can be calculated as
    \begin{equation}
         S_e = \frac{\langle\hat{u} - \tilde{u},\hat{u} - \tilde{u}\rangle_{L^2}}{\langle\hat{u},\hat{u}\rangle_{L^2}}.
    \end{equation}
    For the Euler equations, this indicator was calculated with respect to the density. A threshold value $S_e^*$ was defined as $S_e^* = \epsilon p^{-4}$ for some coefficient $\epsilon$ based on the assumption that, in one dimension, $S_e \sim p^{-4}$ for solutions in $C^0$. A suitable value of $\epsilon$ was empirically determined to be 0.01.  For the RD-FR method, the RD scheme was utilized for elements where $S_e \geqslant S_e^*$, whereas the FR scheme was utilized where $S_e < S_e^*$.

\subsection{Maximum Wavespeed Estimate}
    To calculate $\lambda_{\max}$ in \cref{eq:dij}, a method for approximating the maximum wavespeed for the Riemann problem in \cref{ssec:invariant_sets} is required. In this work, the wavespeed estimate of \citet{Davis1988} is used, calculated as 
    \begin{equation}
        \lambda_{\max} (\mathbf{u}_L, \mathbf{u}_R) = \max \left ( |\mathbf{v}_L\cdot\mathbf{n}| + c_L, |\mathbf{v}_R\cdot\mathbf{n}| + c_R \right ), \quad \mathrm{where} \quad c = \sqrt{\gamma \frac{p}{\rho}}
    \end{equation}
    for the quantities defined in \cref{ex:euler}. In some cases of the Riemann problem, this estimate may not form a true upper bound on the wave speed, and there exist more complex methods for calculating strict or exact upper bounds \citep{Toro2020, Guermond2016b}. However, in the cases to be shown, we found little difference with stricter estimates for the wave speed. 
\section{Results}
\label{sec:results}
    The RD and the RD-FR methods described in \cref{sec:numerical} were applied to a series of standard one and two-dimensional test cases for the Euler equations. In order to more conveniently express the solution for the Euler equations, we define the vector of primitive variables as $\mathbf{q}=[\rho,\mathbf{v},p]^T$ for the conservative variables defined in \cref{ex:euler} with $\gamma = 1.4$.
    
\subsection{Sod Shock Tube}
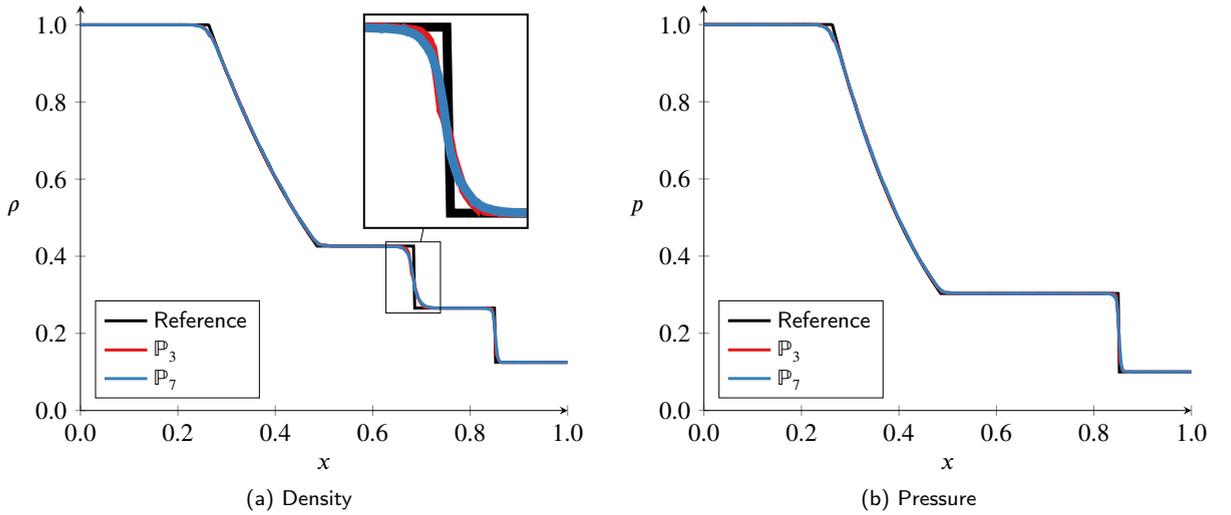
\begin{figure}[tbhp]
        \centering
        \subfloat[Density]{\label{fig:sod_d}\adjustbox{width=0.48\linewidth,valign=b}{    \begin{tikzpicture}[spy using outlines={rectangle, height=3cm,width=2.3cm, magnification=3, connect spies}]
		\begin{axis}[name=plot1,
		    axis line style={latex-latex},
		    axis x line=left,
            axis y line=left,
		    xlabel={$x$},
		    xtick={0,0.2,0.4,0.6,0.8,1},
    		xmin=0,
    		xmax=1,
    		x tick label style={
        		/pgf/number format/.cd,
            	fixed,
            	fixed zerofill,
            	precision=1,
        	    /tikz/.cd},
    		ylabel={$\rho$},
    		ylabel style={rotate=-90},
    		ytick={0,0.2,0.4,0.6,0.8,1},
    		ymin=0,
    		ymax=1.05,
    		y tick label style={
        		/pgf/number format/.cd,
            	fixed,
            	fixed zerofill,
            	precision=1,
        	    /tikz/.cd},
    		legend style={at={(0.03,0.03)},anchor=south west,font=\small},
    		legend cell align={left},
    		style={font=\normalsize}]
    		
			\addplot[color=black, style={very thick}]
				table[x=x,y=r,col sep=comma,unbounded coords=jump]{./figs/data/sod_exact.csv};
    		\addlegendentry{Reference}
    		
			\addplot[color={Set1-A}, style={very thick}]
				table[x=x,y=rho,col sep=comma,unbounded coords=jump]{./figs/data/sod_p3.csv};
    		\addlegendentry{$\mathbb{P}_3$}
    		
			\addplot[color={Set1-B}, style={very thick}]
				table[x=x,y=rho,col sep=comma,unbounded coords=jump]{./figs/data/sod_p7.csv};
    		\addlegendentry{$\mathbb{P}_7$}
			
			\coordinate (spypoint) at (axis cs:0.683,0.345);
            \coordinate (magnifyglass) at (axis cs:0.75,0.75);
		\end{axis} 		

        \spy [black] on (spypoint) in node[fill=white] at (magnifyglass);
		
	\end{tikzpicture}}}
        ~
        \subfloat[Pressure]{\label{fig:sod_p} \adjustbox{width=0.48\linewidth,valign=b}{    \begin{tikzpicture}[spy using outlines={rectangle, height=3cm,width=2.3cm, magnification=3, connect spies}]
		\begin{axis}[name=plot1,
		    axis line style={latex-latex},
		    axis x line=left,
            axis y line=left,
		    xlabel={$x$},
		    xtick={0,0.2,0.4,0.6,0.8,1},
    		xmin=0,
    		xmax=1,
    	    x tick label style={
        		/pgf/number format/.cd,
            	fixed,
            	fixed zerofill,
            	precision=1,
        	    /tikz/.cd},
    		ylabel={$p$},
    		ylabel style={rotate=-90},
    		ytick={0,0.2,0.4,0.6,0.8,1},
    		ymin=0,
    		ymax=1.05,
    		y tick label style={
        		/pgf/number format/.cd,
            	fixed,
            	fixed zerofill,
            	precision=1,
        	    /tikz/.cd},
    		legend style={at={(0.03,0.03)},anchor=south west,font=\small},
    		legend cell align={left},
    		style={font=\normalsize}]
    		
			\addplot[color=black, style={very thick}]
				table[x=x,y=p,col sep=comma,unbounded coords=jump]{./figs/data/sod_exact.csv};
    		\addlegendentry{Reference}
    		
			\addplot[color={Set1-A}, style={very thick}]
				table[x=x,y=p,col sep=comma,unbounded coords=jump]{./figs/data/sod_p3.csv};
    		\addlegendentry{$\mathbb{P}_3$}
    		
			\addplot[color={Set1-B}, style={very thick}]
				table[x=x,y=p,col sep=comma,unbounded coords=jump]{./figs/data/sod_p7.csv};
    		\addlegendentry{$\mathbb{P}_7$}
				
		\end{axis} 		
	\end{tikzpicture}}}
        \caption{\label{fig:sod}Sod shock tube problem at $t=0.2$ with 512 DoF.}
    \end{figure}

    The Sod shock tube problem \citep{Sod1978} is a canonical test case for evaluating the ability of numerical schemes to resolve discontinuities as it contains the three main features of Riemann problems: an expansion fan, contact discontinuity, and shock wave. The problem is solved on the domain $\mathcal{D}=[0,1]$ with the initial condition
    \begin{equation}
        \mathbf{q}(x,0) = \mathbf{q}_0(x) = \begin{cases}
            \mathbf{q}_l, &\mbox{if } x\leqslant 0.5,\\
            \mathbf{q}_r, &\mbox{else},
        \end{cases} \quad \mathrm{given} \quad \mathbf{q}_l = \begin{bmatrix}
            1 \\ 0 \\ 1
        \end{bmatrix}, \quad \mathbf{q}_r = \begin{bmatrix}
            0.125 \\ 0 \\ 0.1
        \end{bmatrix}.
    \end{equation}
    The results of the $\mathbb{P}_3$ and $\mathbb{P}_7$ RD-FR method are shown in \cref{fig:sod} at $t = 0.2$. The degrees of freedom (DoF) were fixed at 512 for both orders, and the reference solution was calculated analytically. The features in the shock tube problem were sufficiently resolved using the RD-FR method without the introduction of oscillations, and no notable degradation in accuracy was observed when extending the scheme to a higher order. It should however be noted that the approach of \citet{Hopf1950} with vanishing viscosity shows that oscillatory solutions for the Euler equations are physical in the presence of dissipation.

    The RD scheme was then applied independently of the FR scheme to show the convergence of the error. The error was calculated with respect to the density, and the $L^1$, $L^2$, and $L^\infty$ norms were defined as
    \begin{align}
        &\|\epsilon\|_{h,1} = \int_\mathcal{D} |\rho - \rho_\mathrm{exact}|\ \mathrm{d} x, \\ 
        &\|\epsilon\|_{h,2} = \sqrt{\int_\mathcal{D} \left(\rho - \rho_\mathrm{exact}\right)^2\ \mathrm{d} x}, \\
        &\|\epsilon\|_{h,\infty} = \max{\left(| \mathbf{\rho} - \mathbf{\rho}_\mathrm{exact} | \right)}.
    \end{align}
    The data of \cref{tab:sod_error_l1,tab:sod_error_l2,tab:sod_error_linf} shows the behavior of the error as the order of the RD scheme was varied in comparison to a first-order finite volume approach, denoted by $\mathbb{P}_0$. In both the $L^1$ and $L^2$ norm, the error converged as expected for a first-order accurate system, but the error for a given resolution decreased with increasing order. This behavior is in contrast to that seen in the method of \citet{Guermond2016} where increasing order causes an increase in error due to the increase in dissipation. An approximately constant maximum error was observed in the $L^\infty$ norm, and the location of the maximum error was within the contact discontinuity which is consistent with the behavior of approximate Riemann solvers for which the contact discontinuity is often the most challenging \citep{Toro1997_4}.

    

    \begin{figure}[tbhp]
        \centering
        \begin{tabular}{r c c c c c c c c}\toprule
	        DoF & $\mathbb{P}_0$ &$\mathbb{P}_1$ &$\mathbb{P}_2$ &$\mathbb{P}_3$ &$\mathbb{P}_4$ &$\mathbb{P}_5$ &$\mathbb{P}_6$ &$\mathbb{P}_7$ \\ \midrule
	        
	        $256$ & \num{1.57E-02} & \num{3.18E-02} & \num{2.46E-02} & \num{1.88E-02} & \num{1.46E-02} & \num{1.24E-02	} & \num{1.04E-02} & \num{8.71E-03} \\ 
	        
	        $512$ & \num{1.03E-02} & \num{2.06E-02} & \num{1.64E-02} & \num{1.24E-02} & \num{9.77E-03} & \num{8.04E-03} & \num{6.72E-03} & \num{5.77E-03} \\
	        
            $1024$ & \num{6.54E-03} & \num{1.31E-02} & \num{1.02E-02} & \num{7.91E-03} & \num{6.24E-03} & \num{5.13E-03} & \num{4.31E-03} & \num{3.78E-03} \\
            
            $2048$ & \num{4.15E-03} & \num{8.35E-03} & \num{6.67E-03} & \num{4.98E-03} & \num{3.90E-03} & \num{3.29E-03} & \num{2.75E-03} & \num{2.40E-03} \\ 
            
            $4096$ & \num{2.64E-03} & \num{5.28E-03} & \num{4.12E-03} & \num{3.17E-03} & \num{2.46E-03} & \num{2.08E-03} & \num{1.73E-03} & \num{1.51E-03} \\
        \textbf{RoC} & $\mathbf{0.6456}$& $\mathbf{0.6484}$& $\mathbf{0.6454}$& $\mathbf{0.6452}$& $\mathbf{0.6463}$& $\mathbf{0.644}$& $\mathbf{0.6465}$& $\mathbf{0.6322}$ \\ \bottomrule
        \end{tabular}
        \captionof{table}{\label{tab:sod_error_l1} $L_1$ norm of the density error for the Sod shock tube problem at varying orders. Rate of convergence shown in bold.}
    \end{figure}
    
    \begin{figure}[tbhp]
        \centering
        \begin{tabular}{r c c c c c c c c}\toprule
	        DoF & $\mathbb{P}_0$ &$\mathbb{P}_1$ &$\mathbb{P}_2$ &$\mathbb{P}_3$ &$\mathbb{P}_4$ &$\mathbb{P}_5$ &$\mathbb{P}_6$ &$\mathbb{P}_7$ \\ \midrule
	        
	        $256$ & \num{2.49E-02} & \num{3.61E-02} & \num{3.23E-02} & \num{2.78E-02} & \num{2.46E-02} & \num{2.28E-02} & \num{2.06E-02} & \num{1.86E-02} \\
	        
	        $512$ & \num{1.89E-02} & \num{2.68E-02} & \num{2.40E-02} & \num{2.10E-02} & \num{1.88E-02} & \num{1.68E-02} & \num{1.57E-02} & \num{1.41E-02} \\
	        
            $1024$ & \num{1.42E-02} & \num{2.00E-02} & \num{1.81E-02} & \num{1.59E-02} & \num{1.40E-02} & \num{1.27E-02} & \num{1.16E-02} & \num{1.10E-02} \\
            
            $2048$ & \num{1.09E-02} & \num{1.56E-02} & \num{1.41E-02} & \num{1.20E-02} & \num{1.08E-02} & \num{9.95E-03} & \num{8.90E-03} & \num{8.42E-03} \\
            
            $4096$ & \num{8.72E-03} & \num{1.23E-02} & \num{1.09E-02} & \num{9.54E-03} & \num{8.42E-03} & \num{7.79E-03} & \num{7.06E-03} & \num{6.51E-03} \\
        \textbf{RoC} & $\mathbf{0.3822}$& $\mathbf{0.3887}$& $\mathbf{0.3902}$& $\mathbf{0.3893}$& $\mathbf{0.3893}$& $\mathbf{0.3854}$& $\mathbf{0.3909}$& $\mathbf{0.3773}$ \\ \bottomrule
        \end{tabular}
        \captionof{table}{\label{tab:sod_error_l2} $L_2$ norm of the density error for the Sod shock tube problem at varying orders. Rate of convergence shown in bold.}
    \end{figure}
    
    \begin{figure}[tbhp]
        \centering
        \begin{tabular}{r c c c c c c c c}\toprule
	        DoF & $\mathbb{P}_0$ &$\mathbb{P}_1$ &$\mathbb{P}_2$ &$\mathbb{P}_3$ &$\mathbb{P}_4$ &$\mathbb{P}_5$ &$\mathbb{P}_6$ &$\mathbb{P}_7$ \\ \midrule
	        
	        $256$ & \num{8.71E-02} & \num{9.25E-02} & \num{9.61E-02} & \num{8.86E-02} & \num{9.19E-02} & \num{8.79E-02} & \num{8.74E-02} & \num{8.46E-02} \\
	        
	        $512$ & \num{8.90E-02} & \num{8.88E-02} & \num{9.12E-02} & \num{9.08E-02} & \num{8.75E-02} & \num{8.93E-02} & \num{9.23E-02} & \num{8.71E-02} \\
	        
            $1024$ & \num{8.79E-02} & \num{8.48E-02} & \num{9.03E-02} & \num{8.99E-02} & \num{8.79E-02} & \num{8.74E-02} & \num{8.68E-02} & \num{8.64E-02} \\
            
            $2048$ & \num{8.69E-02} & \num{9.25E-02} & \num{9.77E-02} & \num{8.53E-02} & \num{8.85E-02} & \num{9.79E-02} & \num{8.70E-02} & \num{8.69E-02} \\
            
            $4096$ & \num{8.60E-02} & \num{8.69E-02} & \num{8.72E-02} & \num{9.01E-02} & \num{8.66E-02} & \num{9.77E-02} & \num{8.72E-02} & \num{8.48E-02} \\
        \textbf{RoC} & -- & -- & -- & -- & -- & -- & -- & -- \\ \bottomrule
        \end{tabular}
        \captionof{table}{\label{tab:sod_error_linf} $L_\infty$ norm of the density error for the Sod shock tube problem at varying orders.}
    \end{figure}

\subsection{Shu--Osher Problem}
\begin{figure}[tbhp]
        \centering
        \subfloat[Density]{\label{fig:shu-osher_density}\adjustbox{width=0.48\linewidth,valign=b}{     \begin{tikzpicture}[spy using outlines={rectangle, height=3cm,width=2.3cm, magnification=3, connect spies}]
		\begin{axis}[name=plot1,
		    axis line style={latex-latex},
		    axis x line=left,
            axis y line=left,
		    xlabel={$x$},
		    xtick={-5,-2.5,0,2.5,5},
    		xmin=-5,
    		xmax=5,
    		x tick label style={
        		/pgf/number format/.cd,
            	fixed,
            	fixed zerofill,
            	precision=1,
        	    /tikz/.cd},
    		ylabel={$\rho$},
    		ylabel style={rotate=-90},
    		ytick={0,1,2,3,4,5},
    		ymin=0,
    		ymax=5,
    		y tick label style={
        		/pgf/number format/.cd,
            	fixed,
            	fixed zerofill,
            	precision=0,
        	    /tikz/.cd},
    		legend style={at={(0.03,0.03)},anchor=south west,font=\small},
    		legend cell align={left},
    		style={font=\normalsize}]
    		
			\addplot[color=black, style={very thick}]
				table[x=x,y=d,col sep=comma,unbounded coords=jump]{./figs/data/osher_p0_2000.csv};
    		\addlegendentry{Reference}
    		
			\addplot[color={Set1-A}, style={very thick}]
				table[x expr={\thisrow{x}-5},y=rho,col sep=comma,unbounded coords=jump]{./figs/data/shu-osher_p3.csv};
    		\addlegendentry{$\mathbb{P}_3$}
    		
			\addplot[color={Set1-B}, style={very thick}]
				table[x expr={\thisrow{x}-5},y=rho,col sep=comma,unbounded coords=jump]{./figs/data/shu-osher_p7.csv};
    		\addlegendentry{$\mathbb{P}_7$}
			
		\end{axis} 		

		
	\end{tikzpicture}}}
        ~
        \subfloat[Pressure]{\label{fig:shu-osher_pressure}\adjustbox{width=0.48\linewidth,valign=b}{    \begin{tikzpicture}[spy using outlines={rectangle, height=3cm,width=2.3cm, magnification=3, connect spies}]
		\begin{axis}[name=plot1,
		    axis line style={latex-latex},
		    axis x line=left,
            axis y line=left,
		    xlabel={$x$},
		    xtick={-5,-2.5,0,2.5,5},
    		xmin=-5,
    		xmax=5,
    		x tick label style={
        		/pgf/number format/.cd,
            	fixed,
            	fixed zerofill,
            	precision=1,
        	    /tikz/.cd},
    		ylabel={$p$},
    		ylabel style={rotate=-90},
    		ytick={0,2,4,6,8,10,12},
    		ymin=0,
    		ymax=12,
    		y tick label style={
        		/pgf/number format/.cd,
            	fixed,
            	fixed zerofill,
            	precision=0,
        	    /tikz/.cd},
    		legend style={at={(0.03,0.03)},anchor=south west,font=\small},
    		legend cell align={left},
    		style={font=\normalsize}]
    		
			\addplot[color=black, style={very thick}]
				table[x=x,y=p,col sep=comma,unbounded coords=jump]{./figs/data/osher_p0_2000.csv};
    		\addlegendentry{Reference}
    		
			\addplot[color={Set1-A}, style={very thick}]
				table[x expr={\thisrow{x}-5},y=p,col sep=comma,unbounded coords=jump]{./figs/data/shu-osher_p3.csv};
    		\addlegendentry{$\mathbb{P}_3$}
    		
			\addplot[color={Set1-B}, style={very thick}]
				table[x expr={\thisrow{x}-5},y=p,col sep=comma,unbounded coords=jump]{./figs/data/shu-osher_p7.csv};
    		\addlegendentry{$\mathbb{P}_7$}
			
		\end{axis} 		

		
	\end{tikzpicture}}}
        \caption{\label{fig:shu-osher}Shu--Osher shock-sine wave interaction problem at $t=0.18$ with 1024 DoF.}
    \end{figure}
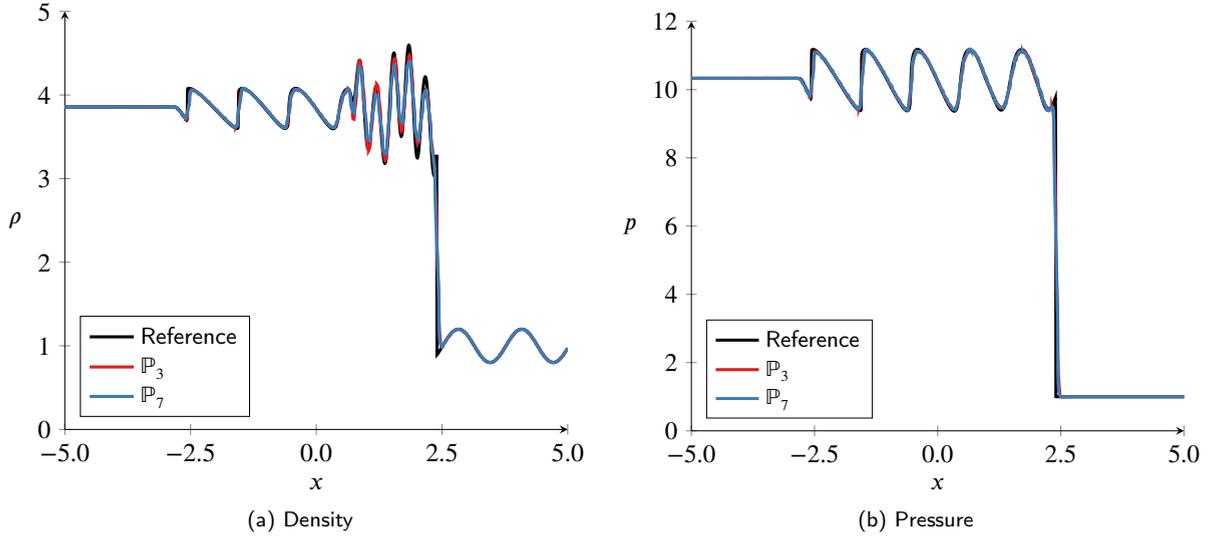

    The case of \citet{Shu1988} tests the ability of the scheme to resolve discontinuities in the presence of physical oscillations. The problem is solved on a domain of $\mathcal{D}=[-5,5]$ with the initial condition
    \begin{equation}
         \mathbf{q}(x,0) = \mathbf{q}_0(x) = \begin{cases}
            \mathbf{q}_l, &\mbox{if } x\leqslant -4, \\
            \mathbf{q}_r, &\mbox{else},
        \end{cases} \quad \mathrm{given} \quad
        \mathbf{q}_l = \begin{bmatrix}
            3.857143 \\ 2.629369 \\ 10.333333
        \end{bmatrix}, \quad
        \mathbf{q}_r = \begin{bmatrix}
            1 + 0.2\sin{5x} \\ 0\\ 1
        \end{bmatrix}.
    \end{equation}
    The density initially contains a sinusoidal oscillation that the shock propagates through. This oscillation can induce instabilities, however, overly dissipative schemes can cause the physical oscillations in the solution to become damped. \cref{fig:shu-osher} shows the results of the $\mathbb{P}_3$ and $\mathbb{P}_7$ RD-FR method with 1024 degrees of freedom at $t = 0.18$. A reference solution was obtained via a highly-resolved exact Godunov-type solver \cite{Toro1997_4}. In both cases, the RD-FR method was able to adequately resolve the discontinuities without excessively dissipating the physical oscillations in the system, although the amplitude of the oscillations was better predicted at a lower order. Furthermore, negligible undershoots were observed in the solution at the discontinuities in both cases.
        
\subsection{Isentropic Euler Vortex}
    The isentropic Euler vortex \citep{Shu1998} is commonly used for verifying the accuracy of a numerical scheme as the results can be directly compared to the analytic solution. The initial conditions are given as 
    \begin{equation}
        \mathbf{q}(\mathbf{x}, 0)  = \begin{bmatrix}
            p^\frac{1}{\gamma} \\
            V_x + \frac{S}{2 \pi R} (y-y_0)\phi(r) \\
            V_y - \frac{S}{2 \pi R} (x-x_0)\phi(r) \\
            \frac{1}{\gamma M^2} \left(1 - \frac{S^2 M^2 (\gamma-1)} {8 \pi^2}\phi(r)^2\right)^\frac{\gamma}{\gamma-1}
        \end{bmatrix}, \quad \mathrm{where} \quad r = \|\mathbf{x}-\mathbf{x}_0\|_2 \quad \mathrm{and} \quad \phi(r) = \exp{\left(\frac{1-r^2}{2R^2}\right)}.
    \end{equation}
    The vortex is characterized by the parameters $S = 13.5$ denoting the strength of the vortex, $R = 1.5$ the radius, and $V_x = 0$, $V_y = 1$ the advection velocities. The freestream Mach number $M$ was set to 0.4, and the boundary conditions were set to periodic for the domain $\mathcal{D} = [-10, 10]^2$ on a uniform quadrilateral mesh. After one convective time in which the vortex has returned to its original position, the $L^2$ norm of the density error was defined as  
    
    \begin{equation}
        \|e\|_2 = \sqrt{\int_\mathcal{D} \left (\rho(x,y) - \rho_0(x,y) \right)^2 \ \mathrm{d} \mathbf{x}}.
    \end{equation}
    
\begin{figure}[tbhp]
        \centering
            \adjustbox{width=0.45\linewidth,valign=b}{\begin{tikzpicture}[spy using outlines={rectangle, height=3cm,width=2.3cm, magnification=3, connect spies}]
	\begin{loglogaxis}[name=plot1,
		xlabel={$1/\overline{h}$},
    	xmin=2e2,xmax=4e3,
    	xtick={1e2,1e3,1e4},
    	ylabel={$\| e \|_2$},
    	ymin=5e-4,ymax=2e-2,
    	grid=both,
    	legend style={at={(0.03,0.03)},anchor=south west,font=\small},
    	legend cell align={left},
    	style={font=\normalsize}]
    	
		\addplot[color=black, style={dotted,very thick},forget plot] coordinates{(4e2,2e-2) (4e3,2e-3)};
		\node [above,color=black] at (axis cs:1.5e3,6e-3) {$\mathcal{O}(\overline{h})$};
		
		\addplot[color={black}, style={very thick}] table[x expr={1./\thisrow{hbar}},y=p0,col sep=comma,unbounded coords=jump]{./figs/data/icv_errors.csv};
		\addlegendentry{$\mathbb{P}_0$};
		
		\addplot[color={Set1-A}, style={very thick}] table[x expr={1./\thisrow{hbar}},y=p1,col sep=comma,unbounded coords=jump]{./figs/data/icv_errors.csv};
		\addlegendentry{$\mathbb{P}_1$ RD};
		
		\addplot[color={Set1-B}, style={very thick}] table[x expr={1./\thisrow{hbar}},y=p2,col sep=comma,unbounded coords=jump]{./figs/data/icv_errors.csv};
		\addlegendentry{$\mathbb{P}_2$ RD};
		
		\addplot[color={Set1-C}, style={very thick}] table[x expr={1./\thisrow{hbar}},y=p3,col sep=comma,unbounded coords=jump]{./figs/data/icv_errors.csv};
		\addlegendentry{$\mathbb{P}_3$ RD};
		
		\addplot[color={Set1-D}, style={very thick}] table[x expr={1./\thisrow{hbar}},y=p4,col sep=comma,unbounded coords=jump]{./figs/data/icv_errors.csv};
		\addlegendentry{$\mathbb{P}_4$ RD};
		
		
	\end{loglogaxis}
\end{tikzpicture}}
        \caption{\label{fig:icv}Density error for the isentropic Euler vortex after one convective time.}
    \end{figure}
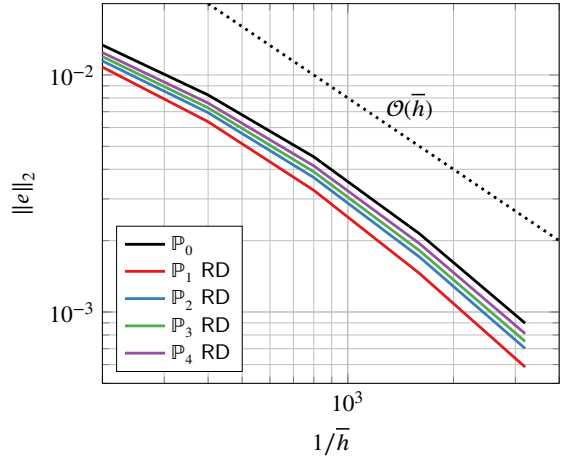

    The density error for a fixed CFL of 0.1 is shown in \cref{fig:icv} with respect to the characteristic element size $\overline{h}$ defined as 
    \begin{equation}
        \overline{h} = \frac{1}{\sqrt{N_{DOF}}},
    \end{equation}
    where $N_{DOF}$ denotes the number of degrees of freedom. The RD scheme at various orders is compared to a first-order finite volume approach, denoted by $\mathbb{P}_0$, at varying values of $\overline{h}$. As expected, the RD scheme is first-order accurate when used independently of a higher-order FR scheme. In comparison to a first-order finite volume approach, the RD scheme offered superior accuracy in terms of the magnitude of the error. The largest benefit in accuracy was attained when using $\mathbb{P}_1$ RD, where up to a 35\% decrease in error was observed, but this benefit was reduced as the order was increased. This was in contrast to the behavior of the error of non-smooth solutions as in the Sod shock tube problem. Although formally first-order accurate, in practice, the RD scheme offered a slightly higher convergence rate than the $\mathbb{P}_0$ approach. This is most evidently observed in the error curve of the $\mathbb{P}_1$ RD method where the relative error reduction monotonically increases from 20\% to 35\% over the sweep of $\overline{h}$.
    
\subsection{Forward Facing Step}
    The forward facing step problem of \citet{Woodward1984} consists of a Mach three flow in a wind tunnel with a step perturbation. The problem is solved on the domain $\mathcal{D} = [0,1] \times [0,3] \setminus [0.6, 3] \times [0, 0.2]$ with the initial condition $\mathbf{q}(x,0) = [1.4, 3, 0, 1]^T$.
    The boundary condition at the inlet was fixed at $\mathbf{q}(x,0)$, while at the outlet, no boundary condition was applied. At the top and bottom boundaries, an adiabatic slip wall condition was enforced such that $\mathbf{n} \cdot \mathbf{v} = 0$. The corner of the step was rounded with a radius of 0.01, and a uniform quadrilateral mesh with a characteristic length $h = 1/200$ was used for the majority of the grid except for the region by the rounded corner where unstructured quadrilaterals of similar size were used. 
    
\begin{figure}[tbhp]
        \centering
        \includegraphics[width=\textwidth]{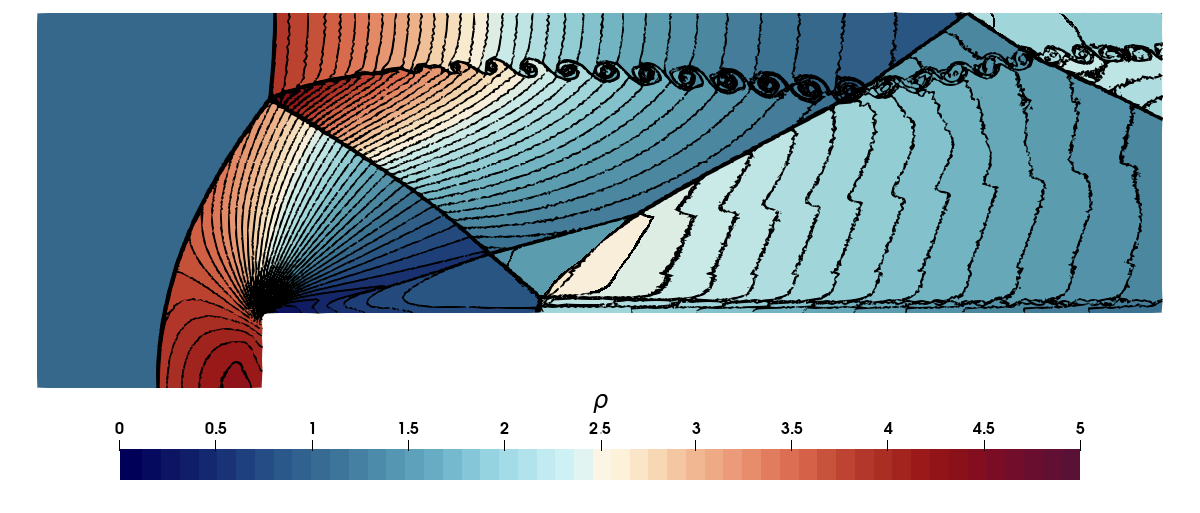}
        \caption{Contours of density for the forward facing step problem at $t = 4$ with $\mathbb{P}_3$ RD-FR and $h = 1/200$.}
        \label{fig:ffs}
    \end{figure}

    The solution of the forward facing step problem at $t = 4$ as predicted by the $\mathbb{P}_3$ RD-FR method is shown in \cref{fig:ffs} using 50 equispaced contours of density. The contours show the interaction of several shock waves as well as the onset of Kelvin-Helmholtz instabilities emanating from the upper shock wave interaction. For more dissipative schemes, these instabilities are difficult to resolve, but the RD-FR scheme was able to predict the rollup of the shear layer and advect the vortices through shock waves without dissipating them. 

\subsection{Richtmyer--Meshkov Instability }
    The Richtmyer--Meskhov instability, predicted analytically by \citet{Richtmyer1960} and shown experimentally by \citet{Meshkov1972}, occurs when a contact discontinuity is acted upon by an impulse, generally as a result of a propagating shock wave. The problem is solved on the domain $\mathcal{D} = [0,10] \times [0,\pi]$ with the initial condition
    \begin{equation}
         \mathbf{q}(x,0) = \mathbf{q}_0(x) = \begin{cases}
            \mathbf{q}_l, &\mbox{if } x\leqslant 1, \\
            \mathbf{q}_c, &\mbox{if } 1 < x \leqslant 3 + a \sin (\omega y),\\
            \mathbf{q}_r, &\mbox{else},
        \end{cases}  \quad \mathrm{given} \quad
        \mathbf{q}_l = \begin{bmatrix}
            1 \\ 0 \\ 0 \\ 1.35
        \end{bmatrix}, \quad
        \mathbf{q}_c = \begin{bmatrix}
            1 \\ 0 \\ 0 \\ 0.1
        \end{bmatrix},\quad
        \mathbf{q}_r = \begin{bmatrix}
            35 \\ 0 \\ 0 \\ 0.1
        \end{bmatrix},
    \end{equation}
    where the parameters $a = 1/4$, $\omega = 4$ dictate the shape of the initial perturbation and the resulting instability behavior \citep{Zanotti2015}. For these initial conditions, the Atwood number was $17/18$, denoting a \textit{light-to-heavy} type Richtmyer--Meshkov problem. The boundary condition was fixed at $\mathbf{q}_l$ at the inlet, periodic at the top and bottom boundaries, and free at the outlet. A uniform quadrilateral mesh with a characteristic length $h =  1/100$ was used. 
    
\begin{figure}[tbhp]
        \centering
        \subfloat[Linear contour map]{\label{fig:rmi_lin}\adjustbox{width=0.45\linewidth,valign=b}{\includegraphics[width=\textwidth]{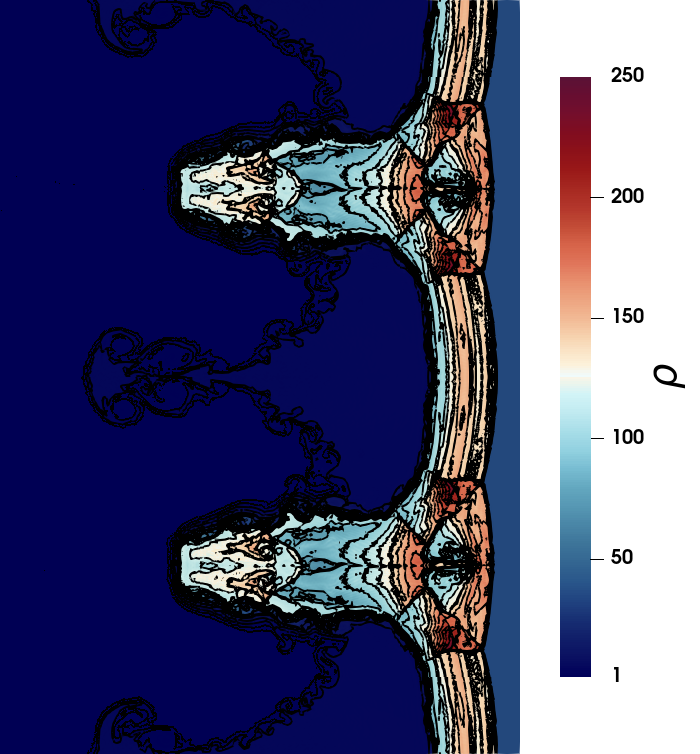}}}
        ~
        \subfloat[Logarithmic contour map]{\label{fig:rmi_log}\adjustbox{width=0.45\linewidth,valign=b}{\includegraphics[width=\textwidth]{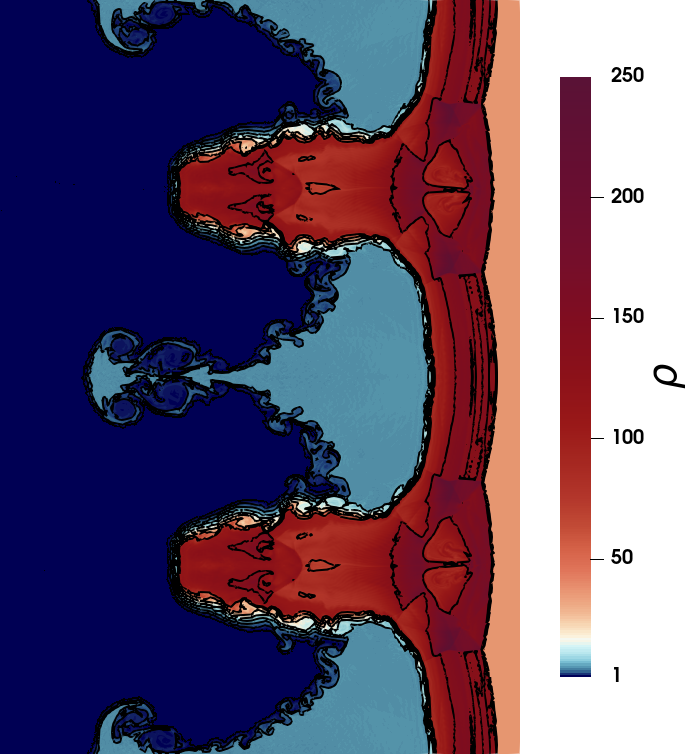}}}
        \caption{\label{fig:rmi}Contours of density for the Richtmyer--Meshkov instability problem at $t = 10$ with $\mathbb{P}_3$ RD-FR and $h = 1/100$.}
    \end{figure}

    The results of the $\mathbb{P}_3$ RD-FR method at $t = 10$ is shown in \cref{fig:rmi} using 20 equispaced and logspaced contours of density. The contour maps show two distinct scales in the problem: shock wave interactions with density jumps of $\mathcal{O}(100)$ and the mushrooming effect of the instability with density jumps of $\mathcal{O}(1)$. After the passage of the incident shock wave, small scale structures of the instability were observed with little smearing of the interface indicating that the numerical diffusion introduced by the scheme does not excessively degrade the accuracy. 
    
\section{Conclusions}
\label{sec:conclusions}
A novel scheme for discontinuous finite element approximations of hyperbolic systems of equations was introduced. The proposed Riemann difference approach uses a staggered grid of nodal solution and flux points within each element, and the flux is calculated by an approximate solution of the Riemann problem posed between the adjacent solution points. For the numerical dissipation introduced by this formulation of the flux polynomial, the scheme was proven to be invariant domain preserving under a presented set of conditions. The method was paired with a higher-order flux reconstruction scheme via a discontinuity sensor and was applied to the Euler equations in one and two dimensions using predominantly affine tensor-product elements.  Numerical experiments showed the ability of the RD scheme to adequately resolve discontinuities without introducing excessive dissipation. When utilized on its own, the scheme was shown to be first-order accurate but with favorable error properties in comparison to first-order finite volume approaches.

\section*{Acknowledgements}
\label{sec:ack}
This research did not receive any specific grant from funding agencies in the public, commercial, or not-for-profit sectors.

\bibliographystyle{unsrtnat}
\bibliography{reference}

\begin{thebibliography}{39}
\providecommand{\natexlab}[1]{#1}
\providecommand{\url}[1]{\texttt{#1}}
\expandafter\ifx\csname urlstyle\endcsname\relax
  \providecommand{\doi}[1]{doi: #1}\else
  \providecommand{\doi}{doi: \begingroup \urlstyle{rm}\Url}\fi

\bibitem[Hopf(1950)]{Hopf1950}
Eberhard Hopf.
\newblock The partial differential equation $u_t + uu_x = \mu_{xx}$.
\newblock \emph{Communications on Pure and Applied Mathematics}, 3\penalty0
  (3):\penalty0 201--230, September 1950.

\bibitem[Godunov(1959)]{Godunov1959}
Sergei~Konstantinovich Godunov.
\newblock A difference method for numerical calculation of discontinuous
  solutions of the equations of hydrodynamics.
\newblock \emph{Matematicheskii Sbornik}, 89\penalty0 (3):\penalty0 271--306,
  1959.

\bibitem[Lax(2006)]{Lax2006}
Peter~D. Lax.
\newblock Gibbs phenomena.
\newblock \emph{Journal of Scientific Computing}, 28\penalty0 (2-3):\penalty0
  445--449, May 2006.

\bibitem[{von Neumann} and Richtmyer(1950)]{VonNeumann1950}
J.~{von Neumann} and R.~D. Richtmyer.
\newblock A method for the numerical calculation of hydrodynamic shocks.
\newblock \emph{Journal of Applied Physics}, 21\penalty0 (3):\penalty0
  232--237, March 1950.

\bibitem[Tadmor(1990)]{Tadmor1990}
Eitan Tadmor.
\newblock Shock capturing by the spectral viscosity method.
\newblock \emph{Computer Methods in Applied Mechanics and Engineering},
  80\penalty0 (1-3):\penalty0 197--208, June 1990.

\bibitem[Glaubitz et~al.(2017)Glaubitz, \"{O}ffner, and Sonar]{Glaubitz2017}
Jan Glaubitz, Philipp \"{O}ffner, and Thomas Sonar.
\newblock Application of modal filtering to a spectral difference method.
\newblock \emph{Mathematics of Computation}, 87\penalty0 (309):\penalty0
  175--207, August 2017.

\bibitem[Persson and Peraire(2006)]{Persson2006}
Per-Olof Persson and Jaime Peraire.
\newblock Sub-cell shock capturing for discontinuous {G}alerkin methods.
\newblock In \emph{44th {AIAA} Aerospace Sciences Meeting and Exhibit}.
  American Institute of Aeronautics and Astronautics, January 2006.

\bibitem[Guermond and Popov(2016{\natexlab{a}})]{Guermond2016}
Jean-Luc Guermond and Bojan Popov.
\newblock Invariant domains and first-order continuous finite element
  approximation for hyperbolic systems.
\newblock \emph{{SIAM} Journal on Numerical Analysis}, 54\penalty0
  (4):\penalty0 2466--2489, January 2016{\natexlab{a}}.

\bibitem[Guermond et~al.(2019)Guermond, Popov, and Tomas]{Guermond2019}
Jean-Luc Guermond, Bojan Popov, and Ignacio Tomas.
\newblock Invariant domain preserving discretization-independent schemes and
  convex limiting for hyperbolic systems.
\newblock \emph{Computer Methods in Applied Mechanics and Engineering},
  347:\penalty0 143--175, April 2019.

\bibitem[Glimm(1965)]{Glimm1965}
James Glimm.
\newblock Solutions in the large for nonlinear hyperbolic systems of equations.
\newblock \emph{Communications on Pure and Applied Mathematics}, 18\penalty0
  (4):\penalty0 697--715, November 1965.

\bibitem[Chueh et~al.(1977)Chueh, Conley, and Smoller]{Chueh1977}
Kai~N Chueh, Charles~C Conley, and Joel~A Smoller.
\newblock Positively invariant regions for systems of nonlinear diffusion
  equations.
\newblock \emph{Indiana University Mathematics Journal}, 26\penalty0
  (2):\penalty0 373--392, 1977.

\bibitem[Hoff(1985)]{Hoff1985}
David Hoff.
\newblock Invariant regions for systems of conservation laws.
\newblock \emph{Transactions of the American Mathematical Society},
  289\penalty0 (2):\penalty0 591--610, 1985.

\bibitem[Boris and Book(1997)]{Boris1997}
Jay~P. Boris and David~L. Book.
\newblock Flux-corrected transport.
\newblock \emph{Journal of Computational Physics}, 135\penalty0 (2):\penalty0
  172--186, August 1997.

\bibitem[Guermond et~al.(2011)Guermond, Pasquetti, and Popov]{Guermond2011}
Jean-Luc Guermond, Richard Pasquetti, and Bojan Popov.
\newblock Entropy viscosity method for nonlinear conservation laws.
\newblock \emph{Journal of Computational Physics}, 230\penalty0 (11):\penalty0
  4248--4267, May 2011.

\bibitem[Kopriva and Kolias(1996)]{Kopriva1996}
David~A. Kopriva and John~H. Kolias.
\newblock A conservative staggered-grid {C}hebyshev multidomain method for
  compressible flows.
\newblock \emph{Journal of Computational Physics}, 125\penalty0 (1):\penalty0
  244--261, April 1996.

\bibitem[Dafermos(2010)]{Dafermos2010_9}
Constantine~M. Dafermos.
\newblock \emph{Hyperbolic Conservation Laws in Continuum Physics}, chapter~9,
  pages 271--324.
\newblock Springer Berlin Heidelberg, 2010.

\bibitem[Lax(1957)]{Lax1957}
P.~D. Lax.
\newblock Hyperbolic systems of conservation laws {II}.
\newblock \emph{Communications on Pure and Applied Mathematics}, 10\penalty0
  (4):\penalty0 537--566, 1957.

\bibitem[Dafermos(2008)]{Dafermos2008}
Constantine Dafermos.
\newblock A variational approach to the {R}iemann problem for hyperbolic
  conservation laws.
\newblock \emph{Discrete and Continuous Dynamical Systems}, 23\penalty0
  (1/2):\penalty0 185--195, September 2008.

\bibitem[Liu et~al.(2004)Liu, Vinokur, and Wang]{Liu2004}
Yen Liu, Marcel Vinokur, and Z.~J. Wang.
\newblock Discontinuous spectral difference method for conservation laws on
  unstructured grids.
\newblock In \emph{Computational Fluid Dynamics}, pages 449--454. Springer
  Berlin Heidelberg, 2004.

\bibitem[Huynh(2007)]{Huynh2007}
H.~T. Huynh.
\newblock A flux reconstruction approach to high-order schemes including
  discontinuous {G}alerkin methods.
\newblock In \emph{18th {AIAA} Computational Fluid Dynamics Conference}.
  American Institute of Aeronautics and Astronautics, June 2007.

\bibitem[Vincent et~al.(2010)Vincent, Castonguay, and Jameson]{Vincent2010}
P.~E. Vincent, P.~Castonguay, and A.~Jameson.
\newblock A new class of high-order energy stable flux reconstruction schemes.
\newblock \emph{Journal of Scientific Computing}, 47\penalty0 (1):\penalty0
  50--72, September 2010.

\bibitem[Witherden et~al.(2014)Witherden, Farrington, and
  Vincent]{Witherden2014}
F.D. Witherden, A.M. Farrington, and P.E. Vincent.
\newblock {PyFR}: An open source framework for solving
  advection{\textendash}diffusion type problems on streaming architectures
  using the flux reconstruction approach.
\newblock \emph{Computer Physics Communications}, 185\penalty0 (11):\penalty0
  3028--3040, November 2014.

\bibitem[Gottlieb et~al.(2001)Gottlieb, Shu, and Tadmor]{Gottlieb2001}
Sigal Gottlieb, Chi-Wang Shu, and Eitan Tadmor.
\newblock Strong stability-preserving high-order time discretization methods.
\newblock \emph{{SIAM} Review}, 43\penalty0 (1):\penalty0 89--112, January
  2001.

\bibitem[Hesthaven and Warburton(2008)]{Hesthaven2008a}
Jan~S. Hesthaven and Tim Warburton.
\newblock \emph{Nodal Discontinuous {G}alerkin Methods}.
\newblock Springer New York, 2008.

\bibitem[Zwanenburg and Nadarajah(2016)]{Zwanenburg2016}
Philip Zwanenburg and Siva Nadarajah.
\newblock Equivalence between the energy stable flux reconstruction and
  filtered discontinuous {G}alerkin schemes.
\newblock \emph{Journal of Computational Physics}, 306:\penalty0 343--369,
  February 2016.

\bibitem[Witherden et~al.(2016)Witherden, Vincent, and Jameson]{Witherden2016}
F.D. Witherden, P.E. Vincent, and A.~Jameson.
\newblock High-order flux reconstruction schemes.
\newblock In \emph{Handbook of Numerical Analysis}, pages 227--263. Elsevier,
  2016.

\bibitem[Rusanov(1962)]{Rusanov1962}
V.V Rusanov.
\newblock The calculation of the interaction of non-stationary shock waves and
  obstacles.
\newblock \emph{{USSR} Computational Mathematics and Mathematical Physics},
  1\penalty0 (2):\penalty0 304--320, January 1962.

\bibitem[Roe(1981)]{Roe1981}
P.L Roe.
\newblock Approximate {R}iemann solvers, parameter vectors, and difference
  schemes.
\newblock \emph{Journal of Computational Physics}, 43\penalty0 (2):\penalty0
  357--372, October 1981.

\bibitem[Davis(1988)]{Davis1988}
S.~F. Davis.
\newblock Simplified second-order {G}odunov-type methods.
\newblock \emph{{SIAM} Journal on Scientific and Statistical Computing},
  9\penalty0 (3):\penalty0 445--473, May 1988.

\bibitem[Toro et~al.(2020)Toro, M\"{u}ller, and Siviglia]{Toro2020}
E.F. Toro, L.O. M\"{u}ller, and A.~Siviglia.
\newblock Bounds for wave speeds in the {R}iemann problem: Direct theoretical
  estimates.
\newblock \emph{Computers {\&} Fluids}, 209:\penalty0 104640, September 2020.

\bibitem[Guermond and Popov(2016{\natexlab{b}})]{Guermond2016b}
Jean-Luc Guermond and Bojan Popov.
\newblock Fast estimation from above of the maximum wave speed in the {R}iemann
  problem for the {E}uler equations.
\newblock \emph{Journal of Computational Physics}, 321:\penalty0 908--926,
  September 2016{\natexlab{b}}.

\bibitem[Sod(1978)]{Sod1978}
Gary~A Sod.
\newblock A survey of several finite difference methods for systems of
  nonlinear hyperbolic conservation laws.
\newblock \emph{Journal of Computational Physics}, 27\penalty0 (1):\penalty0
  1--31, April 1978.

\bibitem[Toro(1997)]{Toro1997_4}
Eleuterio~F. Toro.
\newblock The {R}iemann problem for the {E}uler equations.
\newblock In \emph{{R}iemann Solvers and Numerical Methods for Fluid Dynamics},
  chapter~4, pages 115--157. Springer Berlin Heidelberg, 1997.

\bibitem[Shu and Osher(1988)]{Shu1988}
Chi-Wang Shu and Stanley Osher.
\newblock Efficient implementation of essentially non-oscillatory
  shock-capturing schemes.
\newblock \emph{Journal of Computational Physics}, 77\penalty0 (2):\penalty0
  439--471, August 1988.

\bibitem[Shu(1998)]{Shu1998}
Chi-Wang Shu.
\newblock Essentially non-oscillatory and weighted essentially non-oscillatory
  schemes for hyperbolic conservation laws.
\newblock In \emph{Lecture Notes in Mathematics}, pages 325--432. Springer
  Berlin Heidelberg, 1998.

\bibitem[Woodward and Colella(1984)]{Woodward1984}
Paul Woodward and Phillip Colella.
\newblock The numerical simulation of two-dimensional fluid flow with strong
  shocks.
\newblock \emph{Journal of Computational Physics}, 54\penalty0 (1):\penalty0
  115--173, April 1984.

\bibitem[Richtmyer(1960)]{Richtmyer1960}
Robert~D. Richtmyer.
\newblock Taylor instability in shock acceleration of compressible fluids.
\newblock \emph{Communications on Pure and Applied Mathematics}, 13\penalty0
  (2):\penalty0 297--319, May 1960.

\bibitem[Meshkov(1972)]{Meshkov1972}
E.~E. Meshkov.
\newblock Instability of the interface of two gases accelerated by a shock
  wave.
\newblock \emph{Fluid Dynamics}, 4\penalty0 (5):\penalty0 101--104, 1972.

\bibitem[Zanotti and Dumbser(2015)]{Zanotti2015}
O.~Zanotti and M.~Dumbser.
\newblock High order numerical simulations of the {R}ichtmyer{\textendash}
  {M}eshkov instability in a relativistic fluid.
\newblock \emph{Physics of Fluids}, 27\penalty0 (7):\penalty0 074105, July
  2015.

\end{thebibliography}


\clearpage
\begin{appendices}

\end{appendices}


\end{document}